# GLOBAL EXPONENTIAL STABILIZATION OF ACYCLIC TRAFFIC NETWORKS

Maria Kontorinaki*, Iasson Karafyllis** and Markos Papageorgiou*

*Dynamic Systems and Simulation Laboratory,
Technical University of Crete, Chania, 73100, Greece
(emails: mkontorinaki@dssl.tuc.gr, markos@dssl.tuc.gr)
**Dept. of Mathematics, National Technical University of Athens,
Zografou Campus, 15780, Athens, Greece (email: iasonkar@central.ntua.gr )

**Abstract**

This work is devoted to the construction of explicit feedback control laws for the robust, global, exponential stabilization of general, uncertain, discrete-time, acyclic traffic networks. We consider discrete-time, uncertain network models which satisfy very weak assumptions. The construction of the controllers and the rigorous proof of the robust, global, exponential stability for the closed-loop system are based on recently proposed vector-Lyapunov function criteria, as well as the fact that the network is acyclic. It is shown, in this study, that the latter requirement is necessary for the existence of a robust, global, exponential stabilizer of the desired uncongested equilibrium point of the network. An illustrative example demonstrates the applicability of the obtained results to realistic traffic flow networks.

## 1. Introduction

Networks are large-scale entities representing different types of physical or cyber-physical systems such as fluid flow networks, communication networks, smart grids, etc. [1]–[4]. Particular emphasis is given in this study to traffic networks for which a plethora of diverse infrastructures can be addressed on the basis of a unifying modeling approach (see for example [5]–[7]). More specifically, traffic networks can be modelled as urban road networks consisting of interconnected links which are modelled as store-and-forward components [8] or cell-transmission links [9]; large urban networks consisting of smaller homogeneous sub-networks [10]; freeway networks consisting of series of links, which are modelled, e.g., via general discretized LWR (Lighthill-Whitham-Richards) models [11]–[13] or its simplified CTM (Cell Transmission Model) version [14]; large mixed (corridor) networks consisting of urban and freeway links [15].

Recently, many researchers have addressed the stabilization of equilibrium points of large-scale discrete-time systems. However, the verification of stability for large-scale systems still remains a challenging problem on its own. To this purpose, many tools have been proposed in the literature such as vector Lyapunov functions that are very useful for large-scale discrete-time systems. Sufficient stability conditions by means of vector Lyapunov functions have been proposed in [16] (pages 792-798). In addition, small-gain conditions have been proposed in [17], which can be expressed by means of a vector Lyapunov function formulation (as shown in [18], Chapter 5). Recently in [19], sufficient conditions have been provided for the robust, global, exponential stability of nonlinear, large-scale, uncertain networks by means of vector Lyapunov functions.

The provided results in [19], as it is shown therein, can be easily applied to traffic networks. Traffic networks, satisfying specific assumptions, have also been studied in [5], where sufficient conditions for the local stability of the uncongested equilibrium point are provided; while in [20] the equilibriums of the CTM are analyzed based on monotone systems theory. There are several other works that address stability issues within more specific modeling frameworks for traffic networks. For example, in [21] necessary and sufficient conditions are derived for stable equilibrium accumulations in the undersaturated regimes of macroscopic fundamental diagrams; while [22] studies the stability of equilibriums of a traffic assignment model. However, studies that address rigorously stabilization issues are quite rare. Stability results for simple traffic control systems have been considered in [23] where sufficient conditions for the local and global ISS property of vehicular-traffic networks are provided under the effect of PI regulators. Moreover, in [24] the authors propose link layer feedback (velocity) control laws that stabilize simple, multi-lane and two-dimensional freeway models. Finally, in [13] nonlinear feedback control laws are provided for the robust global exponential stabilization of the uncongested equilibrium point of general freeway models.

In this work, a general model for acyclic networks consisting of an arbitrary number of elementary components with constant turning and exit rates has been developed. The components can be interconnected to form any two-dimensional structure with no cycles for the overall network. Specific instances of the proposed general model result in traffic network structures and problems that can be considered as special cases of the proposed network model and include all the traffic network structures mentioned above. Based on this modeling framework, the results in [19] are utilized for the developed uncertain models of acyclic networks. More specifically, this study provides a parameterized family of explicit feedback control laws which can robustly, globally, exponentially stabilize the desired Uncongested Equilibrium Point (UEP) of a given acyclic traffic network. The achieved stabilization is robust with respect to: i) any uncertainty related to the fundamental diagram of traffic flow; as well as ii) the overall uncertain nature of the developed model when congestion phenomena are present. In fact, in the latter case, the model which describes the time evolution of the network variables is almost completely uncertain (besides the requirement of known and constant turning and exit rates). Furthermore, the assumptions that surround the proposed methodology are weak enough to render the methodology applicable to other kinds of acyclic networks instead of traffic networks. Finally, we emphasize that, as it is proved herein (Proposition 3.1), the requirement regarding the absence of cycles inside the network is utterly necessary for the existence of a robust, global, exponential stabilizer of the UEP of the network.

The provided results generalize some of the results provided in [13] (see Remark 5 and Remark 6 in the following sections) and can be used for ramp-metering control of freeway networks. Preliminary testing and comparison with other existing sophisticated control strategies, provided in [13], demonstrate the efficacy of the proposed methodology as a ramp metering strategy for freeways. However, the present generalized methodology can also be used as perimeter control strategy as well as for arterial (or corridor) networks with arbitrary topology that contain no cycles. To the best of our knowledge, this is the first paper that addresses rigorously global stabilization issues for such problems.

The structure of the present work is as follows. Section 2 includes the model derivation as well as the discussion on the properties and the consequences of the considered modelling framework. The main results of this work are presented in Section 3 while the proofs of the main results can be found in Section 4. An illustrative example of a freeway-to-freeway network is presented in Section 5 and finally the concluding remarks of the paper are given in Section 6.

**Definitions and Notation.** In this paper, we adopt the following notation and terminology:

* $\Re_+ := [0,+\infty)$. $\Re_+^n := (\Re_+)^n$. For every set $S$, $S^n = \underbrace{S \times ... \times S}_{n \text{ times}}$ for every positive integer $n$. For a set $S \subseteq \Re^n$, $\text{int}(S)$ denotes the interior of $S$ (which may be empty).

* By $C^0(A;\Omega)$, we denote the class of continuous functions on $A \subseteq \Re^n$, which take values in $\Omega \subseteq \Re^m$. By $C^k(A;\Omega)$, where $k \geq 1$ is an integer, we denote the class of functions on $A \subseteq \Re^n$ with continuous derivatives of order $k$, which take values in $\Omega \subseteq \Re^m$.

* Let $x, y \in \Re^n$. We say that $x \leq y$ if $(y-x) \in \Re_+^n$ and we say that $x < y$ if $(y-x) \in \text{int}(\Re_+^n)$. The transpose of $x \in \Re^n$ is denoted by $x'$. By $|x|$ we denote the Euclidean norm of $x \in \Re^n$. For every $x \in \Re$, $[x]$ denotes the integer part of $x \in \Re$.

* We denote by $I$ the identity matrix and we denote by $1_{n \times n} \in \Re^{n \times n}$ the matrix for which every entry is equal to one. Moreover, $1_n = (1,...,1)' \in \Re^n$.

* The spectral radius of $\Delta \in \Re^{n \times n}$ is denoted by $\rho(\Delta)$. When all the entries of $\Delta$ are non-negative, then we say that $\Delta$ is non-negative and we write $\Delta \in \Re_+^{n \times n}$.

* We say that the matrix $\Delta \in \Re^{n \times n}$ is upper (lower) triangular if all the entries below (above) the main diagonal are zero. We say that the upper (lower) triangular matrix $\Delta \in \Re^{n \times n}$ is strictly upper (lower) triangular if all the entries of the main diagonal are zero. The diagonal entries of an upper (lower) triangular matrix $\Delta \in \Re^{n \times n}$ are the eigenvalues of $\Delta \in \Re^{n \times n}$.

Let $X \subseteq \Re^n$, $D \subseteq \Re^l$ be non-empty sets and consider the uncertain, discrete-time, dynamical system

$$z^+ = F(d,z), \quad z \in X, \quad d \in D \tag{1.1}$$

where $F: D \times X \to X$ is a mapping. The variable $z \in X$ denotes the state of (1.1) while here (and throughout the paper) $z^+$ denotes the value of the state at the next time instant, i.e., (1.1) expresses the recursive relation $z(t+1) = F(d(t), z(t))$. Let $z^* \in X$ be an equilibrium point of (1.1), i.e., $z^* \in X$ satisfies $z^* = F(d, z^*)$ for all $d \in D$. Notice that the requirement $z^* = F(d, z^*)$ for all $d \in D$ implies that $d \in D$ is a vanishing perturbation, i.e., a disturbance that does not change the position of the equilibrium point of the system.

We use the following definitions throughout the paper.

**Definition 1.1:** *A Trapping Region (TR) for system (1.1) is a set $\Omega \subseteq X$ for which there exists an integer $m \geq 0$ such that for every $z_0 \in X$, $\{d(t) \in D\}_{t=0}^\infty$, the solution $z(t)$ of (1.1) with initial condition $z(0) = z_0$ corresponding to input $\{d(t) \in D\}_{t=0}^\infty$ satisfies $z(t) \in \Omega$ for all $t \geq m$.*

A nonlinear system with a TR is a system for which all solutions enter a specific set after an initial transient period. A direct consequence of Definition 1.1 is that every TR for (1.1) must contain all equilibrium points. We next define the robust, global exponential stability notion for (1.1).

**Definition 1.2:** *We say that $z^* \in X$ is Robustly Globally Exponentially Stable (RGES) for system (1.1) if there exist constants $M, \sigma > 0$ such that for every $z_0 \in X$ and for every sequence $\{d(t) \in D\}_{t=0}^\infty$ the solution $z(t)$ of (1.1) with initial condition $z(0) = z_0$ corresponding to input $\{d(t) \in D\}_{t=0}^\infty$ (i.e., the solution that satisfies $z(t+1) = Z(d(t), z(t))$ for all $t \geq 0$ and $z(0) = z_0$) satisfies the inequality $|z(t) - z^*| \leq M \exp(-\sigma t) |z_0 - z^*|$ for all $t \geq 0$.*

Now, consider the uncertain, discrete-time, control system

$$x^+ = F(x,u), \quad x \in S, \quad u \in U, \tag{1.2}$$

where $F: S \times U \to S$ is a locally bounded mapping and $S \subseteq \Re^n$, $U \subseteq \Re^m$ are non-empty sets. Let $x^* \in S$ be an equilibrium point of (1.2), i.e., there exists $u^* \in U$ so that $x^* = F(x^*, u^*)$. We next define the notion of global asymptotic controllability for (1.2).

**Definition 1.3:** *We say that system (1.2) is globally asymptotically controllable to $x^* \in S$ if for every $x_0 \in S$ there exists $\{u(t) \in U\}_{t=0}^{\infty}$ such that the solution $x(t)$ of (1.2) corresponding to input $\{u(t) \in U\}_{t=0}^{\infty}$ with initial condition $x(0) = x_0$ satisfies $\lim_{t \to \infty} x(t) = x^*$.*

Notice that global asymptotic controllability is a necessary condition for the existence of a globally stabilizing feedback for (1.2) (see [25]).

## 2. Acyclic Networks with Constant Turning and Exit Rates

We consider a generic network which consists of $n$ components (cells). This network may represent a traffic flow network, a fluid flow network or another kind of network. The density of the quantity characterizing each component of the network (e.g. density of vehicles, fluid mass etc.) at time $t \geq 0$ in component $i \in \{1,...,n\}$ is denoted by $x_i(t)$. The outflow and the inflow of the component $i \in \{1,...,n\}$ at time $t \geq 0$ are denoted by $F_{out,i}(t) \geq 0$ and $F_{in,i}(t) \geq 0$, respectively. Consequently, the conservation equation for each component $i \in \{1,...,n\}$ is given by

$$x_i^+ = x_i - F_{out,i} + F_{in,i}, \quad i = 1,...,n, \quad t \geq 0. \tag{2.1}$$

Each component of the network has storage capacity $a_i > 0$ ($i = 1,...,n$). We define

$$S = [0, a_1] \times ... \times [0, a_n], \tag{2.2}$$

which is the set where the state takes values, i.e., $x \in S$. Let $v_i \geq 0$ ($i = 1,...,n$) denote the attempted inflow to component $i \in \{1,...,n\}$ from the region out of the network and set $v = (v_1,...,v_n)' \in \Re_+^n$. Our first assumption is dealing with the outflows. We assume that there exist functions $f_i : D \times [0, a_i] \to \Re_+$, $s_i : D \times S \times \Re_+^n \to [0,1]$ with $f_i(d, x_i) \leq x_i$ for all $(d, x_i) \in D \times [0, a_i]$, where $D \subseteq \Re^l$ is a non-empty, compact set, constants $p_{i,j} \geq 0$, $i, j = 1,...,n$, with $p_{i,i} = 0$ for $i = 1,...,n$, and constants $Q_i \geq 0$, $i = 1,...,n$ so that:

$$\begin{pmatrix} \text{flow from} \\ \text{component } i \text{ to component } j \end{pmatrix} = p_{i,j} s_i(d, x, v) f_i(d, x_i), \text{ for } i, j = 1,...,n, \tag{2.3}$$

$$\begin{pmatrix} \text{flow from} \\ \text{component } i \text{ to regions out of the network} \end{pmatrix} = Q_i s_i(d, x, v) f_i(d, x_i), \text{ for } i = 1,...,n. \tag{2.4}$$

We also assume that:

$$\sum_{j=1}^{n} p_{i,j} + Q_i = 1, \text{ for } i = 1,...,n. \quad (2.5)$$

Some explanations are needed at this point. The functions $f_i : D \times [0, a_i] \to \Re_+$ ($i = 1,...,n$) denote the attempted outflow from the $i$-th cell, i.e., the outflow that will exit the cell if there is sufficient space in the downstream cells. Particularly, the functions $f_i : D \times [0, a_i] \to \Re_+$ ($i = 1,...,n$) remind what in the specialized literature of Traffic Engineering is called the demand-part of the fundamental diagram of the $i$-th cell. In addition, $p_{i,j}$ are turning rates and $Q_i$ are exit rates. The functions $s_i : D \times S \times \Re_+^n \to [0,1]$ ($i = 1,...,n$) are introduced in order to accommodate congestion phenomena. Specifically, these functions assume the value of 1 if the downstream cells can accommodate the whole attempted outflow of the upstream cell; they are less than 1 if the downstream cells cannot accommodate the full attempted outflow, e.g. because they are congested, as it will also explained in more detail later.

Combining (2.3), (2.4), (2.5), we obtain:

$$F_{out,i} = s_i(d, x, v) f_i(d, x_i), \text{ for } i = 1,...,n \quad (2.6)$$

We make the following assumption for the functions $f_i : D \times [0, a_i] \to \Re_+$ ($i = 1,...,n$):

**(H1)** *For each $d \in D$, the function $f_i(d, \cdot) : [0, a_i] \to \Re_+$ satisfies $0 < f_i(d, z) < z$ for all $z \in (0, a_i]$. There exists $\delta_i \in (0, a_i]$ such that for each $d \in D$, the function $f_i(d, \cdot)$ is continuous and increasing on $[0, \delta_i]$. Moreover, there exist constants $L_i \in (0,1)$, $G_i \in (0,1]$, $\tilde{\delta}_i \in (0, \delta_i]$ such that $|f_i(d, z) - f_i(d, y)| \geq L_i |z - y|$ for each $d \in D$ and $y, z \in [0, \tilde{\delta}_i]$ and $|f_i(d, z) - f_i(d, y)| \leq G_i |z - y|$ for each $d \in D$ and $y, z \in [0, \delta_i]$. Finally, there exists a positive constant $f_i^{\min} > 0$ such that for each $d \in D$ it holds that $f_i(d, z) \geq f_i^{\min}$ for all $z \in [\delta_i, a_i]$.*

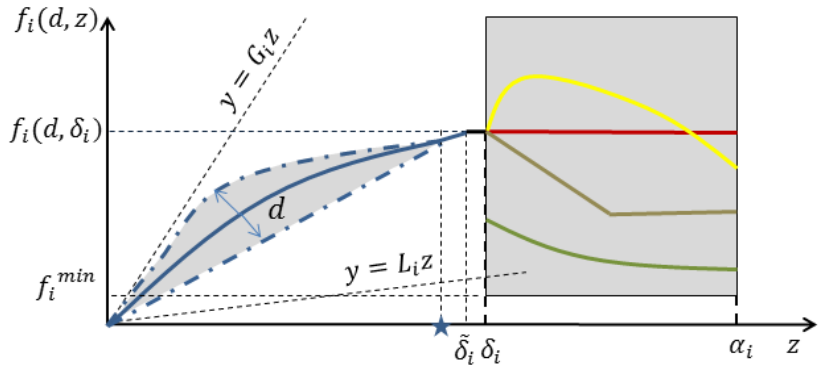

**Figure 1:** Implications of Assumption (H1).

**Remark 1:** Assumption (H1) is a technical assumption that allows a very general class of functions $f_i : D \times [0, a_i] \to \Re_+$ to be taken into account. The implications of assumption (H1) are illustrated in Figure 1. Assumption (H1) includes the basic properties of the so-called "demand function" [12] in the Godunov discretization; $\delta_i$ is the critical density, where $f_i(d, x_i)$ achieves a maximum value (capacity flow). Notice that assumption (H1) includes the possibility of reduced demand flow for overcritical densities (i.e., when $x_i > \delta_i$), since $f_i(d, x_i)$ is allowed to be any arbitrary function (discontinuous or decreasing or, even, increasing), taking any values within the bounds mentioned in (H1) (corresponding to the right grey area in Figure 1), for $x_i > \delta_i$; this could be used to reflect the capacity drop phenomenon of traffic flow, as it is treated in some recent works [26], [27]. Figure 1 presents within the grey area of overcritical densities, four examples of demand functions, which satisfy assumption (H1).

Our second assumption is dealing with the inflows. We assume that there exist functions $g_i \in C^0(D \times S; \Re_+)$, $w_i : D \times S \times \Re_+^n \to [0,1]$ ($i=1,...,n$) with $0 < g_i(d,x) \le a_i - x_i$ for all $(d,x) \in D \times S$ with $x_i < a_i$ and $i=1,...,n$, so that

$$F_{in,i} = w_i(d,x,v)v_i + \sum_{j=1}^{n} p_{j,i} s_j(d,x,v) f_j(d,x_j) \le g_i(d,x), \text{ for all } i=1,...,n \text{ and } (d,x,v) \in D \times S \times \Re_+^n. \quad (2.7)$$

If $v_i + \sum_{j=1}^{n} p_{j,i} f_j(d,x_j) \le g_i(d,x)$, for all $i=1,...,n$ then $w_i(d,x,v) = s_i(d,x,v) = 1$, for $i=1,...,n$. (2.8)

Again, for the case of traffic networks, the functions $g_i : D \times S \to \Re_+$ ($i=1,...,n$) remind what in the specialized literature of Traffic Engineering is called the supply function of the $i$-th cell. When $w_i(d,x,v) + \sum_{j=1}^{n} p_{j,i} s_j(d,x,v) < 1 + \sum_{j=1}^{n} p_{j,i}$ then we say that the $i$-th cell is *congested*. The functions $w_i : D \times S \times \Re_+^n \to [0,1]$ and $s_i : D \times S \times \Re_+^n \to [0,1]$ ($i=1,...,n$) are introduced so that for each cell: (i) the demand is always less than the supply (this is inequality (2.7)), and (ii) when the maximum value of all demands can be accommodated then no congestion phenomena are present (this is implication (2.8)). Priority rules for each junction can be expressed by means of the functions $w_i : D \times S \times \Re_+^n \to [0,1]$ and $s_i : D \times S \times \Re_+^n \to [0,1]$ ($i=1,...,n$).

For traffic flow networks, the supply function is usually given by the function $g_i(d,x) = \min(q_i, c_i(a_i - x_i))$, where $q_i$ represents the maximum admissible inflow of the $i$-th cell and $c_i \in (0,1]$ represents the normalized congestion wave speed. Then, the fundamental diagram (FD) of cell $i$ is composed by the increasing function $f_i(d,x_i)$ for $x_i \in [0,\delta_i]$ and by the non-increasing function $g_i(d,x) = \min(q_i, c_i(a_i - x_i))$ for $x_i \in [\delta_i, a_i]$. Notice here that the uncertainty $d \in D$ has been introduced in order to accommodate the uncertain nature of the fundamental diagram.

Combining equations (2.1), (2.6) and (2.7) we obtain the following uncertain discrete-time system:

$$x_i^+ = x_i + w_i(d,x,v)v_i - s_i(d,x,v)f_i(d,x_i) + \sum_{j=1}^{n} p_{j,i} s_j(d,x,v) f_j(d,x_j), \text{ for } i=1,...,n. \quad (2.9)$$

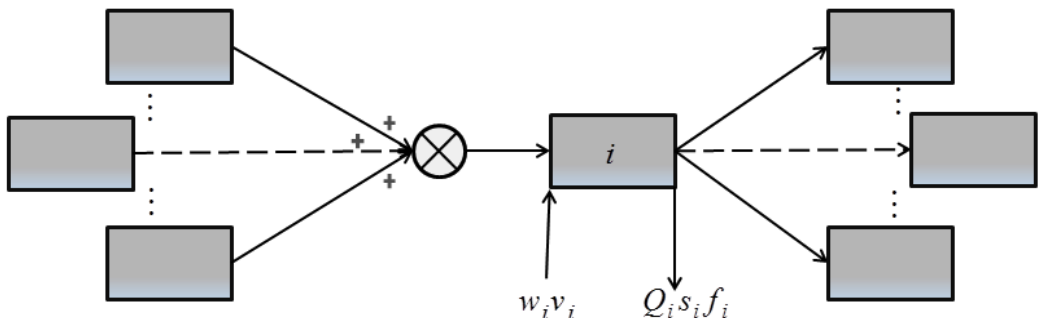

**Figure 2:** Scheme of the network model (2.9).

Figure 2 illustrates schematically the network described by the model (2.9). For physical reasons, we would expect a network of the form (2.9) under assumption (H1) to satisfy the following three properties:

1) If the attempted external inflows $v_i \geq 0$ ($i = 1,...,n$) are "small" for a sufficiently large time period then the network densities will eventually be "small".
2) If $x_i \neq 0$ for some $i = 1,...,n$, then there is at least one non-zero outflow.
3) If the attempted external inflows $v_i \geq 0$ ($i = 1,...,n$) and the $x_i \geq 0$ ($i = 1,...,n$) are "small", then no congestion phenomena are present in the network.

Indeed, consider a network with zero external inflows. If the network does not satisfy property 1 above then it is possible that the network retains a certain amount of density (i.e., the vehicles do not exit). The same situation would occur in the case where property 2 above does not hold. Of course, there are "special" cases (e.g. a gridlock around a cycle) where vehicles are trapped in the network and do not exit, but it is clear that in such situations one cannot deal with congestion phenomena via inflow control, i.e. by making the external inflows sufficiently small. Property 3 is another empirical fact that should be verified to enable inflow control: congestion phenomena are present only when the attempted external inflows $v_i \geq 0$ ($i = 1,...,n$) and the network densities $x_i \geq 0$ ($i = 1,...,n$) are "sufficiently large".

In the aim of guaranteeing that the considered network models actually possess the above properties, we consider only acyclic networks via the following assumption.

**(H2)** *The matrix $P = \{p_{i,j} : i, j = 1,...,n\} \in [0,1]^{n \times n}$ which contains the turning rates of the acyclic network (2.9) is strictly upper triangular.*

**Remark 2:** From a graph-theoretic point of view, directed acyclic graphs are graphs whose vertices can admit a topological sorting. This means, that their vertices can be ordered in such a way, that the starting endpoint of every edge (joining two vertices) occurs earlier in the ordering than the ending endpoint of the edge. Assigning the vertices of the graph to the components or cells of the network, for any given acyclic network, and by using the previous definition, we are in a position to reorder the cells of the network into a topological sorting. The main consequence of this sorting is that the matrix $P = \{p_{i,j} : i, j = 1,...,n\} \in [0,1]^{n \times n}$ containing the turning rates of the network becomes strictly upper triangular [28]–[30].

The following technical lemmas are useful for the analysis of the networks. Their proofs are provided in the Appendix.

**Lemma 2.1:** *For every non-negative, strictly upper triangular matrix $P$ with $\sum_{j=1}^{n} p_{i,j} \leq 1$ for all $i = 1,...,n$, there exist positive constants $r_i > 0$ ($i = 1,...,n$), such that*

$$r_i > \sum_{j=1}^{n} r_j p_{i,j}, \text{ for every } i = 1,...,n. \tag{2.10}$$

**Lemma 2.2:** *Let $L_i \in (0,1)$ and $G_i \in (0,1]$ with $L_i \leq G_i$ for $i = 1,...,n$ be constants and let $P$ be a non-negative, strictly upper triangular matrix with $\sum_{j=1}^{n} p_{i,j} \leq 1$ for $i = 1,...,n$. Then, there exist constants $\xi_i > 0$ ($i = 1,...,n$), such that*

$$\sum_{j=1}^{n} p_{j,i} G_j \xi_j < L_i \xi_i, \text{ for every } i = 1,...,n. \tag{2.11}$$

Using vector notation, inequality (2.11) becomes

$$(P'diag(G) - diag(L))\xi \leq 0$$

which is also equivalent to

$$(I + P'diag(G) - diag(L))\xi \leq \xi. \quad (2.12)$$

**Lemma 2.3:** *Let $L_i \in (0,1)$ and $G_i \in (0,1]$ with $L_i \leq G_i$ for $i = 1,...,n$ be constants and let $P$ be a non-negative, strictly upper triangular matrix with $\sum_{j=1}^{n} p_{i,j} \leq 1$ for $i = 1,...,n$. Then the matrix $I + P'diag(G) - diag(L)$ is a lower triangular matrix with $\rho(I + P'diag(G) - diag(L)) < 1$, where $G = (G_1,...,G_n)$ and $L = (L_1,...,L_n)$.*

The following assumption is a technical assumption, which is related to Property 2 above.

**(H3)** *There exist functions $\tilde{s}_i \in C^0(D \times S \times \Re_+^n;[0,1])$ with $s_i(d,x,v) \geq \tilde{s}_i(d,x,v)$ for all $(d,x,v) \in D \times S \times \Re_+^n$, $i = 1,...,n$, and constants $v_i^{\max} > 0$ ($i = 1,...,n$) such that the following implication holds:*

$$\text{If } x_i \tilde{s}_i(d,x,v) = 0, i = 1,...,n \text{ and } v_i < v_i^{\max}, i = 1,...,n \text{ then } x = 0. \quad (2.13)$$

**Remark 3:** Assumption (H3) guarantees that the functions $s_i : D \times S \times \Re_+^n \to [0,1]$, for $i = 1,...,n$, which have been introduced in model (2.9) in order to accommodate congestion phenomena, should admit a continuous and positive definite lower bound for some $i = 1,...,n$. Implication (2.13) guarantees that if the outflow of every cell of the network is zero then the density of every cell should be zero (Property 2).

We next show that assumption (H3) in conjunction with assumption (H1) and (H2) guarantees that the network (2.9) satisfies Properties 1, 2 above.

**Proposition 2.4:** *Consider the network (2.9) under assumptions (H1), (H2), (H3). Then for every constants $r_i > 0$ ($i = 1,...,n$) satisfying (2.10) and for every family of constants $\tilde{\varepsilon}_i \in \left(0, \min\left(v_i^{\max}, \min\{g_i(d,0) : d \in D\}\right)\right)$ ($i = 1,...,n$), there exists a constant $C > 0$ such that*

$$\left(\sum_{i=1}^{n} r_i x_i\right)^+ \leq (1-C)\sum_{i=1}^{n} r_i x_i + \sum_{i=1}^{n} r_i v_i \quad (2.14)$$

*for all $(d,x) \in D \times S$ and for all $v_i \geq 0$ with $v_i \leq \min\left(v_i^{\max}, \min\{g_i(d,0) : d \in D\}\right) - \tilde{\varepsilon}_i$ ($i = 1,...,n$).*

Inequality (2.14) and induction allows us to show that for every $\omega > 0$ and for sufficiently small external inflows ($v_i(t) \geq 0$ for $i = 1,...,n$ with $v_i(t) \leq \min\left(v_i^{\max}, \min\{g_i(d,0) : d \in D\}\right) - \tilde{\varepsilon}_i$ for all $t \geq 0$) there exists $T > 0$ sufficiently large such that the following estimate holds for all $t \geq T$ for the solution of (2.9), for every initial condition $x(0) \in S$ and for every input $\{d(t) \in D\}_{t=0}^{\infty}$:

$$\sum_{i=1}^{n} r_i x_i(t) \leq \omega + C^{-1} \max_{i=1,...,n} \left(\sup\{v_i(t) : t \geq 0\}\right) \sum_{i=1}^{n} r_i.$$

The above inequality shows that if the attempted external inflows $v_i \geq 0$ ($i = 1,...,n$) are "small" for a sufficiently large time period then the network densities will eventually be "small". This is Property 1 stated above. Property 2 above is a direct consequence of (2.6), (2.13) and the fact that $f_i(d,x_i) = 0 \Leftrightarrow x_i = 0$ (a consequence of Assumption (H1)). Property 3 is a direct consequence of the following assumption and (2.8).

**(H4)** *There exist constants* $\mu_i \in (0, \tilde{\delta}_i)$, $v_i^{\max} > 0$ ($i = 1,...,n$), *such that*

$$v_i^{\max} + \sum_{j=1}^{n} p_{j,i} f_j(d, x_j) \leq g_i(d, x) \text{ for all } i = 1,...,n, \ (d, x) \in D \times S \text{ with } x \leq \mu \quad (2.15)$$

*where* $\mu = (\mu_1,..., \mu_n)'$.

**Remark 4:** Assumption (H4) is a reasonable assumption: if the network densities are small (below a critical value, here denoted by $\mu_i$) and the attempted external inflows are small (below a given $v_i^{\max}$), then the total attempted inflow should be accommodated by the $i$-th cell.

Assumptions (H1), (H2), (H3) and (H4) have important consequences; some of them have been already discussed while the rest are presented in the next section. Those assumptions may fit to many kinds of networks of the form (2.9). In particular, for freeway traffic flow networks the aforementioned assumptions are relatively mild. In fact, the following remark shows that the above assumptions are indeed satisfied for the general freeway models in [13].

**Remark 5:** The above assumptions are satisfied for the general freeway models in [13] due to the following facts:

- A similar assumption as (H1) is used in [13], but assumption (H1) is more general, since uncertain demand functions and relaxed assumptions both for the left, $x_i \leq \delta_i$, and the right, $x_i > \delta_i$, branch of the demand functions have been considered within assumption (H1) in this study.

- The entries of the matrix $P$ which contains the turning rates of the network (2.9) for the freeway models in [13] are as follows

$$\begin{aligned} p_{i,j} &= 0 \text{ for every } j \neq i+1 \text{ with } i = 1,...,n \text{ and } j = 1,...,n \\ p_{i,j} &> 0 \text{ for } j = i+1 \text{ with } i = 1,...,n-1 \text{ and } j = 1,...,n \end{aligned} \quad (2.16)$$

    thus we conclude that assumption (H2) holds for (2.9) with (2.16).

- In [13] the supply functions are given as $g_i(d, x_i) = \min\{q_i, c_i(a_i - x_i)\}$, where $q_i > 0$, $c_i \in (0,1]$ ($i = 1,...,n$) are constants. It can be shown that assumption (H4) is satisfied with $v_i^{\max} > 0$, $\mu_i \in (0, \tilde{\delta}_i)$, $i = 1,...,n$, given by the following recursive formulas

$$\mu_n = \tilde{\delta}_n / 2, \ \mu_i = \min\left(\frac{\tilde{\delta}_i}{2}, \frac{1}{2} \min(q_{i+1}, c_{i+1}(a_{i+1} - \mu_{i+1}))\right) (i = 1,...,n-1) \text{ and} \quad (2.17)$$

$$v_i^{\max} = \frac{1}{2} \min(q_i, c_i(a_i - \mu_i)) \text{ for } i = 1,...,n. \quad (2.18)$$

- It can be also shown that assumption (H3) holds with $v_i^{\max} > 0$ ($i = 1,...,n$) given by (2.18) and the continuous functions $\tilde{s}_i(d, x, v)$ ($i = 1,...,n$) given by

$$\tilde{s}_i(d, x, v) := \min\left(1, \frac{\max(0, \min\{q_{i+1}, c_{i+1}(a_{i+1} - x_{i+1})\} - v_{i+1})}{p_{i,i+1} a_i}\right) (i = 1,...,n-1) \text{ and } \tilde{s}_n := 1. \quad (2.19)$$

## 3. Main Results

Consider a network of the form (2.9) under assumptions (H1), (H2), (H3), (H4). We next assume the existence of a point $x^* = (x_1^*,..., x_n^*)' \in S$ with $x_i^* \in (0, \mu_i)$ for $i = 1,..., n$ and a vector $v^* = (v_1^*,..., v_n^*)' \in \Re_+^n$ with $v_i^* < v_i^{max}$ for $i = 1,..., n$, that satisfy the following equations:

$$f_i(d, x_i^*) = v_i^* + \sum_{j=1}^{n} p_{j,i} f_j(d, x_j^*), \text{ for all } i = 1,..., n \text{ and } d \in D. \quad (3.1)$$

Since $x_i^* \in (0, \mu_i)$, $v_i^* < v_i^{max}$ for $i = 1,..., n$, it follows from (2.15) that the following inequalities hold:

$$v_i^* + \sum_{j=1}^{n} p_{j,i} f_j(d, x_j^*) < g_i(d, x^*), \text{ for all } i = 1,..., n \text{ and } d \in D. \quad (3.2)$$

The point $x^* = (x_1^*,..., x_n^*)' \in S$ is called the UEP (uncongested equilibrium point) of the network corresponding to the vector of external inflows $v^* = (v_1^*,..., v_n^*)' \in \Re_+^n$. Notice that the input $d \in D$ is a vanishing perturbation for system (2.9) with $v(t) \equiv v^*$. This is also illustrated in Figure 1, which shows that the input $d \in D$ does not change the position of the equilibrium point (denoted by a star).

One of the most important consequences of the existence of an UEP and assumptions (H1), (H2), (H3), (H4) is presented below. More specifically, the following proposition reveals the reason for studying acyclic networks (explicitly guaranteed by (H2)) and shows that if the network contains cycles, then the system is not globally asymptotically controllable to the UEP. That means that assumption (H2) is utterly necessary in order to proceed to the study of the stabilization of the network (2.9) because otherwise there is no feedback control law which can render the UEP globally exponentially stable. (Note that proofs of the main results are provided in Section 4 and the Appendix.)

**Proposition 3.1:** *Consider the network (2.9) under assumptions (H1), (H3), (H4). Assume the existence of a point $x^* = (x_1^*,..., x_n^*)' \in S$ with $x_i^* \in (0, \mu_i)$ for $i = 1,..., n$ and a vector $v^* = (v_1^*,..., v_n^*)' \in \Re_+^n$ with $v_i^* < \min(v_i^{max}, \min\{g_i(d,0): d \in D\})$ for $i = 1,..., n$, that satisfy equations (3.1). Assume that the network contains at least one cycle. Then, system (2.9) with input $(v,d) \in \Re_+^n \times D$ is not globally asymptotically controllable to the equilibrium point $x^* = (x_1^*,..., x_n^*)' \in S$.*

We next assume that some of the external inflows may be controlled. Let $b \in \Re_+^n$ be a vector with $b \leq v^*$, let $K \in \Re_+^{n \times n}$ be a non-negative, constant matrix and let $\tau > 0$ be a constant. We set:

$$v = v^* - diag(v^* - b)\left(1_n - h\left(1_n - \tau^{-1} K h(x - x^*)\right)\right) \quad (3.3)$$

where $h: \Re^n \to \Re_+^n$ is the mapping defined by

$$h(x) = \left(\max(0, x_1),..., \max(0, x_n)\right)' \in \Re_+^n, \text{ for all } x \in \Re^n. \quad (3.4)$$

Notice that if $b_i = v_i^*$ for some $i \in \{1, 2,..., n\}$ then it follows from (3.3) that $v_i = v_i^*$, i.e., the external inflow $v_i$ is uncontrolled. Therefore, by assuming (3.3), we have taken into account all possible cases for the control of external inflows. We intend to prove the following theorem, which shows that the UEP can be robustly, globally, exponentially stabilized by the continuous feedback law (3.3), which regulates certain or all the external inflows.

**Theorem 3.2:** *Consider the network (2.9) under assumptions (H1), (H2), (H3), (H4). Assume the existence of a point $x^* = (x_1^*,...,x_n^*)' \in S$ with $x_i^* \in (0, \mu_i)$ for $i = 1,...,n$ and a vector $v^* = (v_1^*,...,v_n^*)' \in \Re_+^n$ with $v_i^* < \min\left(v_i^{\max}, \min\{g_i(d,0) : d \in D\}\right)$ for $i = 1,...,n$, that satisfy equations (3.1). Then there exists an index set $R \subseteq \{i \in \{1,...,n\} : v_i^* > 0\}$, a matrix $K \in \Re_+^{n \times n}$ and a vector $b \in \Re_+^n$ with $0 < b_i < v_i^*$ for $i \in R$, $b_i = v_i^*$ for $i \notin R$ such that for every $\tau \in (0,1)$, $x^* = (x_1^*,...,x_n^*)' \in S$ is Robustly Globally Exponentially Stable for the closed-loop system (2.9) with (3.3).*

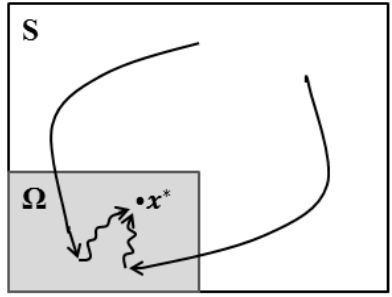

**Figure 3:** Idea behind Theorem 3.2.

Theorem 3.2 is an existence result. However, its proof is constructive and provides formulae (or sufficient conditions) for all constants and for the index set $R$ (see Section 4 and the Appendix). Notice that the index set $R$ is the set of all inflows that must be controlled in order to be able to guarantee that the UEP is RGES. The importance of Theorem 3.2 lies on the following facts:
  a) It provides a family of robust, global, exponential stabilizers (parameterized by $\tau \in (0,1)$) and an explicit feedback law (formula (3.3)).
  b) The achieved stabilization is robust with respect to:
      i. The uncertain nature (introduced by $d \in D$) of the fundamental diagram of traffic flow (by considering uncertain demand and supply functions, $f_i(d,\cdot)$ and $g_i(d,\cdot)$, respectively).
      ii. The overall uncertain nature of the model (2.9) when congestion phenomena are present (by considering uncertain functions $s_i(d,\cdot,\cdot)$ and $w_i(d,\cdot,\cdot)$, with respect to $d \in D$).

Notice here, that the only requirements regarding the functions $s_i(d,\cdot,\cdot)$ (and $w_i(d,\cdot,\cdot)$ respectively) are summarized within the implication (2.8) and assumption (H3). However, implication (2.8) is not a strict requirement since it allows the functions $s_i(d,\cdot,\cdot)$ (for $i = 1,...,n$) to take any value within $[0,1]$, when at least one cell is congested. One possibility for the uncertainty within the functions $s_i(d,\cdot,\cdot)$ (for $i = 1,...,n$) is to be represented with respect to unknown and even time-varying priority rules in the junctions of the network as in [13], where freeway models are considered (which are special cases of (2.9)); however, here, this type of uncertainty may be enhanced by considering priority rules for all the internal inflows of the network. Notice also, that the only requirements regarding the functions $f_i(d,\cdot)$ and $g_i(d,\cdot)$ are summarized within assumption (H1) and the inequality $g_i(d,x) \leq a_i - x_i$ which again allow for a large variety of fundamental diagrams to be considered (see the illustrative example in Section 5).

The main idea behind the proof of Theorem 3.2 is the construction of a vector Lyapunov function for the closed-loop system. The construction of the vector Lyapunov function is based on the existence of a TR $\Omega$ for the system (2.9) in which no congestion phenomena are present. The appropriate selection of the gain matrix $K \in \Re_+^{n \times n}$ in (3.3) forces the selected control action to lead the state in the set $\Omega$ (see Figure 3). In other words, the control action will first eliminate all congestion phenomena and then will drive the state to the desired equilibrium.

**Remark 6:** It is important also to notice that Theorem 3.2 is a generalization of the corresponding theorem in [13] (Theorem 2.1), which shows that a continuous, robust, global, exponential stabilizer exists for the aforementioned freeway models proposed therein. We have already shown (Remark 5) that the considered assumptions in this study are generalizations of the corresponding assumptions and definitions in [13]. But also the feedback stabilizer defined by (3.3) is generalization of the feedback law proposed in [13]. This can be shown by selecting the matrix $K \in \Re_+^{n \times n}$ as $K_{i,j} = \sigma^j$ for every $i, j = 1,...,n$ where $\sigma \in (0,1]$ is a parameter. Therefore, clearly the present study is a comprehensive generalization of [13], and the proposed results can be directly applied to freeway models.

The following proposition shows the existence of a positively invariant region for (2.9).

**Proposition 3.3:** *Consider the network (2.9) under assumptions (H1), (H2), (H3), (H4). Assume the existence of a point $x^* = (x_1^*,..., x_n^*)' \in S$ with $x_i^* \in (0, \mu_i)$ for $i = 1,...,n$ and a vector $v^* = (v_1^*,..., v_n^*)' \in \Re_+^n$ with $v_i^* < v_i^{\max}$ for $i = 1,...,n$, that satisfy equations (3.1). Then there exist constants $\beta_i \in (x_i^*, \mu_i]$ ($i = 1,...,n$) such that for every $b \in \Re_+^n$ with $b \le v^*$, $K \in \Re_+^{n \times n}$ and, $\tau > 0$ it holds that*

$$x \in \Omega, d \in D \Rightarrow x^+ \in \Omega, \quad (3.5)$$

*where $\Omega = [0, \beta_1] \times ... \times [0, \beta_n]$, $h: \Re^n \to \Re_+^n$ is the mapping defined by (3.4) and $x^+$ is given by (2.9) with (3.3).*

Implication (3.5) shows that $\Omega \subset S$ is a positively invariant region for inputs that satisfy $d(t) \in D$ and $0 \le v(t) \le v^* - diag(v^* - b)\left(1_n - h\left(1_n - \tau^{-1} Kh(x(t) - x^*)\right)\right)$ for all $t \ge 0$. It should be noticed that $x^* \in int(\Omega)$, i.e., the UEP is in the interior of the positively invariant region. In order to study the stability properties of the UEP of the network (2.9), we need the following technical lemmas. Their proofs are provided in the Appendix.

**Lemma 3.4:** *Consider the network (2.9) under assumptions (H1), (H2), (H3), (H4). Assume the existence of a point $x^* = (x_1^*,..., x_n^*)' \in S$ with $x_i^* \in (0, \mu_i)$ for $i = 1,...,n$ and a vector $v^* = (v_1^*,..., v_n^*)' \in \Re_+^n$ with $v_i^* < v_i^{\max}$ for $i = 1,...,n$, that satisfy equations (3.1). Then there exist constants $\beta_i \in (x_i^*, \mu_i]$ ($i = 1,...,n$) such that for every $b \in \Re_+^n$ with $b \le v^*$, $K \in \Re_+^{n \times n}$ and $\tau > 0$, implication (3.5) holds and such that*

$$x \in \Omega, d \in D \Rightarrow h(x^+ - x^*) \le \left(I + P'diag(G) - diag(L)\right)h(x - x^*) \quad (3.6)$$

$$x \in \Omega, d \in D \Rightarrow h(x^* - x^+) \le \left(I + P'diag(G) - diag(L)\right)h(x^* - x) + diag(v^* - b)\tau^{-1}Kh(x - x^*) \quad (3.7)$$

*where $\Omega = [0, \beta_1] \times ... \times [0, \beta_n]$, $h: \Re^n \to \Re_+^n$ is the mapping defined by (3.4), $L = (L_1,..., L_n)' \in \Re^n$, $G = (G_1,..., G_n)' \in \Re^n$, $P \in \Re^{n \times n}$ is the matrix $P = \{p_{i,j} : i, j = 1,...,n\}$ and $x^+$ is given by (2.9) with (3.3).*

**Lemma 3.5:** *Consider the network (2.9) under assumptions (H1), (H2), (H3), (H4). Assume the existence of a point $x^* = (x_1^*,..., x_n^*)' \in S$ with $x_i^* \in (0, \mu_i)$ for $i = 1,...,n$ and a vector $v^* = (v_1^*,..., v_n^*)' \in \Re_+^n$ with $v_i^* < v_i^{\max}$ for $i = 1,...,n$, that satisfy equations (3.1). Then there exist constants $\beta_i \in (x_i^*, \mu_i]$ ($i = 1,...,n$) such that for every $b \in \Re_+^n$ with $b \le v^*$, $K \in \Re_+^{n \times n}$ and $\tau > 0$, implications (3.5), (3.6), (3.7) hold and there exists a constant $M > 0$ (depending on $b \in \Re_+^n$, $K \in \Re_+^{n \times n}$ and $\tau > 0$), which satisfies the following property*

$$x \in S, d \in D \Rightarrow |x^+ - x^*| \le M|x - x^*|, \quad (3.8)$$

*where $x^+$ is given by (2.9) with (3.3).*

The following lemma shows the existence of a TR for system (2.9).

**Lemma 3.6:** *Consider the network (2.9) under assumptions (H1), (H2), (H3), (H4). Assume the existence of a point $x^* = (x_1^*,...,x_n^*)' \in S$ with $x_i^* \in (0, \mu_i)$ for $i = 1,...,n$ and a vector $v^* = (v_1^*,...,v_n^*)' \in \Re_+^n$ with $v_i^* < \min\left(v_i^{\max}, \min\{g_i(d,0): d \in D\}\right)$ for $i = 1,...,n$, that satisfy equations (3.1). Let $r = (r_1,...,r_n)' \in \text{int}(\Re_+^n)$ be a vector of constants satisfying (2.10) and let $C > 0$ be the corresponding constant for which inequality (2.14) holds for all $(d,x) \in D \times S$ and for all $v_i \geq 0$ with $v_i \leq v_i^*$ ($i = 1,...,n$). Assume that there exist $b \in \Re_+^n$ with $b \leq v^*$ such that*

$$r'b \leq C \min_{i=1,...,n}(r_i x_i^*). \tag{3.9}$$

*Then there exist constants $\beta_i \in (x_i^*, \mu_i]$ ($i = 1,...,n$) and a matrix $K \in \Re_+^{n \times n}$ such that for every $\tau \in (0,1)$ implications (3.5), (3.6), (3.7) hold and the set $\Omega = [0, \beta_1] \times ... \times [0, \beta_n]$ is a TR for the closed-loop system (2.9) with (3.3).*

Summarizing the above results, the following theorem shows that the UEP is robustly, globally, exponentially stable for the system (2.9) under the proposed feedback regulator (3.3). The proof of Theorem 3.7 is based on the construction of a vector Lyapunov function [19]. Theorem 3.7 is utilized in order to prove the main result of this section, i.e., Theorem 3.2.

**Theorem 3.7:** *Consider the network (2.9) under assumptions (H1), (H2), (H3), (H4). Assume the existence of a point $x^* = (x_1^*,...,x_n^*)' \in S$ with $x_i^* \in (0, \mu_i)$ for $i = 1,...,n$ and a vector $v^* = (v_1^*,...,v_n^*)' \in \Re_+^n$ with $v_i^* < \min\left(v_i^{\max}, \min\{g_i(d,0): d \in D\}\right)$ for $i = 1,...,n$, that satisfy equations (3.1). Let $r = (r_1,...,r_n)' \in \text{int}(\Re_+^n)$ be a vector of constants satisfying (2.10) and let $C > 0$ be the corresponding constant for which inequality (2.14) holds for all $(d,x) \in D \times S$ and for all $v_i \geq 0$ with $v_i \leq v_i^*$ ($i = 1,...,n$). Assume that there exist $b \in \Re_+^n$ with $b \leq v^*$ such that (3.9) holds. Then there exist constants $\beta_i \in (x_i^*, \mu_i]$ ($i = 1,...,n$) and a matrix $K \in \Re_+^{n \times n}$ such that for every $\tau \in (0,1)$, implications (3.5), (3.6), (3.7) hold and the equilibrium point $x^* = (x_1^*,...,x_n^*)' \in S$ is Robustly Globally Exponentially Stable for the closed-loop system (2.9) with (3.3).*

# 4. Proof of Main Results

## 4.1. Proof of Proposition 3.1

Let the index set $E \subseteq \{1,...,n\}$ be the set of all the indices of the cells that are in one of the cycles in the network (2.9). Let also $e \leq n$ be the cardinality of the set $E$. Then, we define $E := \{i_1, i_2,..., i_e\}$ so that $p_{i_1,i_2}, p_{i_2,i_3},..., p_{i_{e-1},i_e}, p_{i_e,i_1} \neq 0$. Moreover, consider an initial condition $x(0)$ for which $x_{i_k}(0) = a_{i_k}$ for every $k = 1,...,e$ (but otherwise arbitrary) and let $\{d(t) \in D\}_{t=0}^{\infty}$ and $\{v(t) \in U\}_{t=0}^{\infty}$ be arbitrary sequences. Due to the fact that $g_{i_k}(d,x) = 0$ if $x_{i_k} = a_{i_k}$ (direct consequence of continuity of $g_i(d,x)$ and the fact that $0 < g_i(d,x) \leq a_i - x_i$ for every $i = 1,...,n$) for every $k = 1,...,e$, we have from (2.7) that for $k = 1$

$$F_{in,i_1}(0) = 0 \Rightarrow w_{i_1}(d(0),x(0),v(0))v_{i_1}(0) + \sum_{j=1}^{n} p_{j,i_1} s_j(d(0),x(0),v(0)) f_j(d(0),x_j(0)) = 0 \Rightarrow w_{i_1}(d(0),x(0),v(0)) = 0$$

and $\sum_{j=1}^{n} p_{j,i_1} s_j(d(0),x(0),v(0)) f_j(d(0),x_j(0)) = 0$.

But the fact that $\sum_{j=1}^{n} p_{j,i_1} s_j(d(0),x(0),v(0)) f_j(d(0),x_j(0)) = 0$ implies that:

$$\begin{aligned} p_{i_e,i_1} s_{i_e}(d(0),x(0),v(0)) f_{i_e}(d(0),x_{i_{e-1}}(0)) &= p_{i_e,i_1} s_{i_e}(d(0),x(0),v(0)) f_{i_e}(d(0),a_{i_e}) \\ &= p_{i_e,i_1} s_{i_e}(d(0),x(0),v(0)) f_{i_e}^{\min} = 0 \Rightarrow s_{i_e}(d(0),x(0),v(0)) = 0. \end{aligned}$$

Repeating the above process for every $k = 1,...,e$, we obtain that $w_{i_k}(d(0),x(0),v(0)) = 0$ and $s_{i_k}(d(0),x(0),v(0)) = 0$ for every $k = 1,...,e$. Therefore, we conclude from (2.9) that $x_{i_k}(1) = x_{i_k}(0) = a_{i_k}$. Using induction, it follows that $x_{i_k}(t) = a_{i_k}$ for every $t \geq 0$. The proof is complete.  ◁

## 4.2. Proof of Proposition 3.3

Lemma 2.2 guarantees that there exists $\xi \in \text{int}(\Re_+^n)$ such that $\sum_{j=1}^{n} p_{j,i} G_j \xi_j < L_i \xi_i$, for $i = 1,...,n$. Using (2.12) and the fact that $h(x) \leq \xi$ for all $x \in \Re^n$ with $x \leq \xi$, we have

$$(I + P' \text{diag}(G) - \text{diag}(L))h(x) < \xi. \tag{4.1}$$

Since $x_i^* \in (0, \mu_i)$ for $i = 1,...,n$, there exists a constant $\varepsilon^* > 0$, sufficiently small, such that $x^* + \varepsilon^* \xi \leq \mu$, where $\mu = (\mu_1,...,\mu_n)' \in \text{int}(\Re_+^n)$. We define:

$$\beta := x^* + \varepsilon^* \xi. \tag{4.2}$$

Let arbitrary $x \in \Omega, d \in D, v \in \Re_+^n$ with $v \leq v^* - \text{diag}(v^* - b)(1_n - h(1_n - \tau^{-1} Kh(x - x^*)))$ be given. Since $x_i \leq \beta_i \leq \mu_i$ and $v_i \leq v_i^* < v_i^{\max}$ for $i = 1,...,n$, it follows from (2.8) and (2.15) that

$$x_i^+ = x_i + v_i - f_i(d, x_i) + \sum_{j=1}^{n} p_{j,i} f_j(d, x_j), \text{ for } i = 1,...,n. \tag{4.3}$$

Using the fact that $v_i \leq v_i^* - (v_i^* - b_i)\left(1 - \max\left(0, 1 - \tau^{-1} \sum_{j=1}^{n} K_{i,j} \max(0, x_j - x_j^*)\right)\right)$ for $i = 1,...,n$, in conjunction with (3.1), we obtain from (4.3):

$$\begin{aligned} x_i^+ &\leq x_i - (v_i^* - b_i) + (v_i^* - b_i) \max\left(0, 1 - \tau^{-1} \sum_{j=1}^{n} K_{i,j} \max(0, x_j - x_j^*)\right) \\ &+ f_i(d, x_i^*) - f_i(d, x_i) + \sum_{j=1}^{n} p_{j,i}\left(f_j(d, x_j) - f_j(d, x_j^*)\right) \end{aligned}, \text{ for } i = 1,...,n. \tag{4.4}$$

Using Assumption (H1), we get $L_i(x_i - x_i^*) \leq f_i(d, x_i) - f_i(d, x_i^*) \leq G_i(x_i - x_i^*)$ for $i = 1,...,n$ and $x_i \geq x_i^*$. Notice that assumption (H1) and the fact that $\mu_i \in (0, \tilde{\delta}_i)$ ($i = 1,...,n$), guarantee that the mappings $x_i \to x_i - f_i(d, x_i)$, $x_i \to f_i(d, x_i)$ are non-decreasing for $x_i \in [0, \beta_i]$, $i = 1,...,n$. It follows that:

$$\begin{aligned} f_i(d, x_i) - f_i(d, x_i^*) &\leq G_i \max(0, x_i - x_i^*), \quad x_i - x_i^* + f_i(d, x_i^*) - f_i(d, x_i) \leq (1 - L_i) \max(0, x_i - x_i^*) \\ &\text{for } x_i \in [0, \beta_i], \, i = 1,...,n. \end{aligned} \tag{4.5}$$

Combining (4.4), (4.5) we obtain for $i = 1,...,n$:

$$x_i^+ \leq x_i^* - (v_i^* - b_i) + (v_i^* - b_i)\max\left(0, 1 - \tau^{-1}\sum_{j=1}^n K_{i,j}\max(0, x_j - x_j^*)\right)$$
$$+ (1 - L_i)\max(0, x_i - x_i^*) + \sum_{j=1}^n p_{j,i}G_j\max(0, x_j - x_j^*) \quad (4.6)$$

Using vector notation and definition (3.4), we are in a position to write inequalities (4.6) in the following form:

$$x^+ \leq x^* - v^* + b + diag(v^* - b)h(1_n - \tau^{-1}Kh(x - x^*)) + (I + P'diag(G) - diag(L))h(x - x^*). \quad (4.7)$$

In order to show (3.5), it suffices to show that

$$x^* - v^* + b + diag(v^* - b)h(1_n - \tau^{-1}Kh(x - x^*)) + (I + P'diag(G) - diag(L))h(x - x^*) \leq \beta,$$
$$\text{for all } x \in \Re_+^n \text{ with } x \leq \beta$$

or equivalently, using (4.2),

$$diag(v^* - b)h(1_n - \tau^{-1}Kh(x - x^*)) + (I + P'diag(G) - diag(L))h(x - x^*) \leq \varepsilon^* \xi + v^* - b,$$
$$\text{for all } x \in \Re_+^n \text{ with } x \leq x^* + \varepsilon^* \xi. \quad (4.8)$$

Setting $x = x^* + \varepsilon^* \zeta$, where $\zeta \in \Re^n$, and using the fact that $h(\varepsilon^* \zeta) = \varepsilon^* h(\zeta)$ for all $\zeta \in \Re^n$ (a direct consequence of definition (3.4)), it follows that (4.8) holds provided that

$$diag(v^* - b)h(1_n - \varepsilon^*\tau^{-1}Kh(\zeta)) + \varepsilon^*(I + P'diag(G) - diag(L))h(\zeta) \leq \varepsilon^* \xi + v^* - b,$$
$$\text{for all } \zeta \in \Re^n \text{ with } \zeta \leq \xi. \quad (4.9)$$

However, inequality (4.1) and the fact that $h(1_n - \varepsilon^*\tau^{-1}Kh(\zeta)) \leq 1_n$ imply (4.9).
The proof is complete.    ◁

### 4.3. Proof of Theorem 3.7

A direct application of Theorem 2.3 in [19]. Indeed, Lemma 3.5 and Lemma 3.6 guarantee that all assumptions of Theorem 2.3 in the above paper hold for the closed-loop system (2.9), (3.3) with

$$V_i(x) := \max(0, x_i - x_i^*) \text{ for } i = 1,...,n \text{ and } V_i(x) := \max(0, x_i^* - x_i) \text{ for } i = n+1,...,2n. \quad (4.10)$$

Notice that definitions (4.10) guarantee the inequality

$$\frac{1}{\sqrt{n}}|x - x^*| \leq \max_{i=1,...,2n}(V_i(x)) = \max_{i=1,...,n}(|x_i - x_i^*|) \leq |x - x^*|, \text{ for all } x \in S \quad (4.11)$$

while inequalities (3.6), (3.7) and definitions (4.10) imply the inequality

$$V(x^+) \leq \Gamma V(x), \text{ for all } (d, x) \in D \times \Omega, \quad (4.12)$$

where $V(x) = (V_1(x),...,V_{2n}(x))' \Re^{2n}$ and

$$\Gamma := \begin{bmatrix} I + P'diag(G) - diag(L) & 0 \\ diag(v^* - b)\tau^{-1}K & I + P'diag(G) - diag(L) \end{bmatrix}. \tag{4.13}$$

Lemma 2.3 guarantees that the matrix $I + P'diag(G) - diag(L)$ is a lower triangular matrix with $\rho(I + P'diag(G) - diag(L)) < 1$. Then, it follows that the matrix $\Gamma$, as defined by (4.13), is a lower triangular matrix with its diagonal entries being the same with the diagonal entries of the matrix $I + P'diag(G) - diag(L)$. Therefore, $\rho(\Gamma) < 1$. The proof is complete. ◁

## *4.4. Proof of Theorem 3.2*

Without loss of generality, by virtue of Theorem 3.7, it suffices to show the existence of $b \in \Re_+^n$ with $b \leq v^*$ such that (3.9) hold. We set:

$$R := \{i \in \{1,...,n\} : v_i^* > 0\} \tag{4.14}$$

$$b := \lambda v^* \tag{4.15}$$

$$\lambda := \min\left(\frac{1}{2}, \frac{C \min_{i=1,...,n}(r_i x_i^*)}{r'v^*}\right). \tag{4.16}$$

Notice that definitions (4.14), (4.15), (4.16) guarantee that (3.9) holds.
The proof is complete. ◁

## 5. Illustrative Example

Consider a 3-lane freeway-to-freeway traffic network of the form (2.9) with $n = 8$ cells. The traffic network consists of two smaller freeways, 2 km each; the first is composed by the cells $i = 1,2,3,4$, and the second is composed by the cells $i = 5,6,7,8$ (see, Figure 4). The cells are homogeneous, each cell being 0.5 km in length. The whole network admits two external inflows; one external inflow at the upstream boundary of the first cell and one external inflow at the upstream boundary of the fifth cell, while there are no intermediate external inflows ($v_2 = v_3 = v_4 = v_6 = v_7 = v_8 = 0$ and $v_1, v_5 \neq 0$). At the end of the first freeway (4th cell) there is an off-ramp joining the second freeway which becomes an on-ramp for the second freeway at the upstream boundary of the 7th cell (see Figure 4). According to this configuration, the exit and turning rates of the freeway are defined as follows

$$Q_i = 0 \text{ for } i = \{1,...8\} \setminus \{4,8\}, \ Q_4 = 0.5, \ Q_8 = 1 \text{ and}$$
$$p_{i,j} = \begin{cases} 1 & \text{if } j = i+1 \text{ and } i \in \{1,...,7\} \setminus \{4\} \\ 0.5 & \text{if } i = 4 \text{ and } j = 7 \\ 0 & \text{if otherwise} \end{cases}, \ i,j = 1,...,8. \tag{5.1}$$

Consequently, the only control possibilities are the inflows $v_1, v_5$. It should be noted here that the 7th cell is a bottleneck for the overall network due to the ramp that joins both freeways. Congestion may be created in the 7th cell, due to high on-ramp demand from the 1st and the 5th cells, and spill back to both freeways depending on the priority rules.

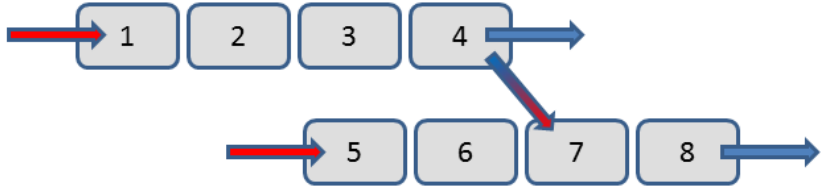

**Figure 4:** The scheme of the freeway-to-freeway network.

All the following simulation tests have been conducted using the following form of the model (2.9), which is expressed by means of the supply function $g_i$ ($i = 1,...,n$):

$$x_i^+ = x_i - s_i(d,x,v)f_i(d,x,v) + \min\left(g_i(d,x), v_i + \sum_{j=1}^{n} p_{j,i} f_j(d,x_j)\right) \quad (5.2)$$

$$s_i(d,x,v) = \begin{cases} \min\left(1, \max\left(0, \dfrac{g_{i+1}(d,x) - v_{i+1}}{p_{i,i+1} f_i(d,x_i)}\right)\right) & \text{if } x_i > 0 \\ 1 & \text{if } x_i = 0 \end{cases}, \text{ for } i \neq 4,8\,,$$

$$s_4(d,x,v) = \begin{cases} \min\left(1, \max\left(0, \dfrac{g_7(d,x) - p_{6,7} f_6(d,x_6)}{p_{4,7} f_4(d,x_4)}\right)\right) & \text{if } x_4 > 0 \\ 1 & \text{if } x_4 = 0 \end{cases} \text{ and } s_8(d,x,v) = 1. \quad (5.3)$$

Notice that, according to (5.3), constant priority rules for the junctions have been taken into account by assuming a full priority rate for the external inflows and by assuming that the mainstream flow coming from the 6$^{th}$ cell has full priority over the mainstream flow coming from the 4$^{th}$ cell. Furthermore, we assume that the simulation time step is $T = 15$ s. However, since all flows and densities are measured in [veh] (as imposed by the form of the model (2.9) and (5.2), (5.3)), the cell length, the time step and the number of lanes do not appear explicitly, but they are only reflected implicitly in the values of every variable and every constant (e.g. critical density, jam density, flow capacity, wave speed etc.) corresponding to density or flow. Appropriate transformations in common traffic units are given for the most critical variables wherever it is needed.

The demand and the supply functions have been defined so as to reflect the uncertainty, $d$, derived from the fundamental diagram of traffic flow. More specifically, we assume that the demand functions are given as a convex combination of several functions $\phi_i$ (e.g., linear or quadratic) satisfying assumption (H1). Furthermore, it should be noted that the functions $\phi_i$ should guarantee that the uncertainty $d$ is a vanishing perturbation for the system (5.2), (5.3), i.e., it does not change the position of the UEP (see Figure 5). Here, six different functions are used to represent the uncertainty in the demand functions. Specifically, the functions $\phi_i$, for $i = 1,...,6$ are given by

$$\phi_1(z) = \frac{5}{11}z, \quad \phi_2(z) = -(13.5/3025)z^2 + 0.7z, \quad \phi_3(z) = (14/3025)z^2 + 0.2z,$$

$$\phi_4(z) = \begin{cases} (-49/3025)z^2 + 0.9z & \text{if } z \in [0, 27.5] \\ (-38/3025)z^2 + 82/55\,z - 19 & \text{if } z \in (27.5, 55] \end{cases},$$

$$\phi_5(z) = \begin{cases} (7/756.25)z^2 + 0.2z & \text{if } z \in [0, 27.5] \\ (21/6050)z^2 + (71.5/1210)z + 8.25 & \text{if } z \in (27.5, 55] \end{cases}, \quad (5.4)$$

$$\phi_6(z) = -\frac{3}{23}z + \frac{740}{23}, \quad \phi_7(z) = (83/52900)z^2 - (4471/10580)z + 46019/1058.$$

Then, the demand functions are given by

$$f_i(d, x_i) = \begin{cases} d_1\phi_1(x_i) + d_2(1-d_1)\phi_2(x_i) + (1-d_2)(1-d_1)\phi_3(x_i) & \text{if } x_i \in [0, 55+2\varepsilon] \\ d_3\phi_6(x_i) + (1-d_3)\phi_7(x_i) & \text{if } x_i \in (55+2\varepsilon, 170] \end{cases}, \text{ for} \quad (5.5)$$

$$i = 1,2,3,4,7,8,$$

$$f_i(d, x_i) = \begin{cases} d_1\phi_1(x_i) + d_2(1-d_1)\phi_4(x_i) + (1-d_2)(1-d_1)\phi_5(x_i) & \text{if } x_i \in [0, 55+2\varepsilon] \\ d_3\phi_6(x_i) + (1-d_3)\phi_7(x_i) & \text{if } x_i \in (55+2\varepsilon, 170] \end{cases}, \text{ for } i = 5,6, \quad (5.6)$$

where $d_i \in [0,1]$, for $i = 1,2,3$, correspond to time-varying weight parameters and $\varepsilon = 10^{-5}$. According to (5.5), (5.6), each cell has the same critical density $\delta_i = 55 + 2\varepsilon$ [veh] ($i=1,...,8$) (corresponding to 36.7 [veh/km/lane] with the above settings) and the same jam density $a_i = 170$ [veh] ($i=1,...,8$) (corresponding to 113.3 [veh/km/lane]). Notice also that, according to (5.5), (5.6), decreasing functions have been considered for overcritical densities, as proposed in [27], so as to incorporate into the model (5.2), (5.3) the capacity drop phenomenon.

As it has already been mentioned, for traffic flow networks the supply functions are usually described by the functions $g_i(d,x) = \min(q_i, c_i(a_i - x_i))$, where $q_i > 0$ represents the maximum inflow for the $i^{th}$ cell and $c_i \in (0,1]$ represents the normalized congestion wave speed.

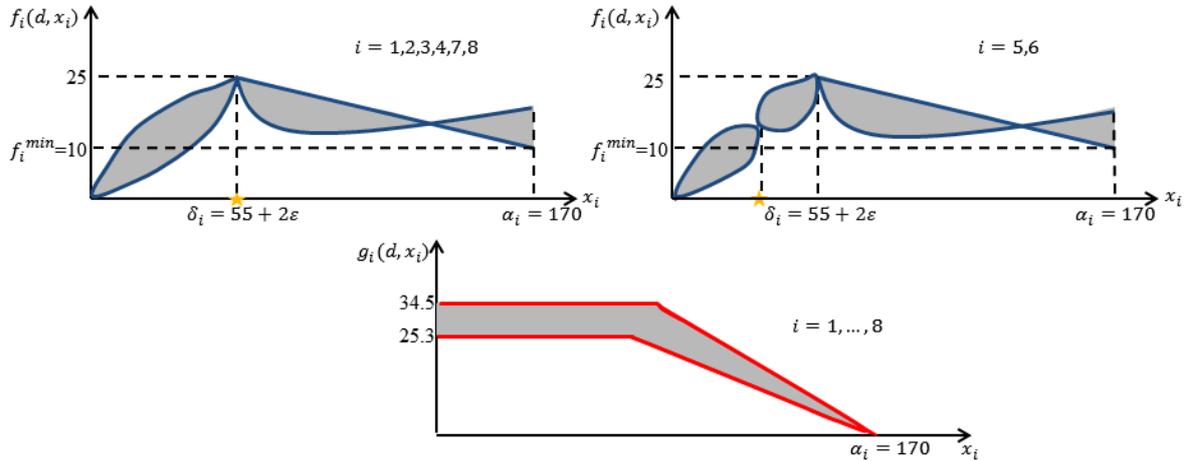

**Figure 5:** Specification of the parameters of the demand and the supply functions of every cell.

Here, in order to consider the uncertainty of the supply functions, we assume that

$$g_i(d, x_i) = d_4 \min(115, a_i - x_i), \quad (5.7)$$

where $d_4 \in [0.22, 0.30]$ is a time-varying parameter resulting to a congestion wave speed within approximately 26 to 36 [km/h] and a maximum inflow approximately between 2000 to 2750 [veh\h\lane]. For the overall system (5.2), (5.3), the uncertainty $d(t) = (d_1(t),...,d_4(t)) \in D$ is a time-varying parameter taking values from a uniform distribution within $D = [0,1]^3 \times [0.22, 0.3]$. Figure 5 visualizes a grey area which includes any possible demand and supply functions.

Assumption (H1) is satisfied for $\delta_i = \tilde{\delta}_i = 55 + 2\varepsilon$ [veh], $f_i^{\min} = 10$ [veh] ($i=1,...,8$), $L_i = 0.2$, $G_i = 0.71$ for $i = 1,2,3,4,7,8$ and $L_i = 0.009$, $G_i = 0.9$ for $i = 5,6$. The cell flow capacities are approximately $f_i(d, \delta_i) = 25$ [veh] ($i = 1,...,8$), corresponding to 2000 [veh/h/lane]. Notice also that the matrix $P$, given by (5.1), satisfies assumption (H2). Assumption (H4) holds for $\mu_i = 55 + \varepsilon$ for $i \neq 5,6$ and $\mu_i = 27.5 + \varepsilon$ for $i = 5,6$, where $\varepsilon = 10^{-5}$ and $v_i^{max} = 0.3$ for $i \neq 1,5$ and $v_1^{max} = v_5^{max} = 25$. Finally, assumption (H3) holds for

$$\tilde{s}_i(d,x,v) := \min\left(1, \frac{\max(0, g_{i+1}(d,x) - v_{i+1})}{p_{i,i+1}a_i}\right) \text{ for } i \neq 4,8,$$

$$\tilde{s}_4(d,x,v) := \min\left(1, \frac{\max(0, g_7(d,x) - v_7 - p_{6,7}f_6(d,x_6))}{p_{4,7}a_4}\right) \text{ and } \tilde{s}_8(d,x,v) := 1 \tag{5.8}$$

and $v_i^{max}$ (for $i = 1,...,8$) as previous.

Here, $R = \{1,5\}$ and therefore we select $b_1 = b_5 = 0.5$ while $b_i = 0$ for every $i \neq 1,5$. Our goal is to globally exponentially stabilize the system at an UEP which is as close as possible to the critical density (due to the fact that the flow value at the critical density is the largest). Equation (3.1) and inequality (3.2) are satisfied by selecting $v^* = (25,0,0,0,12.5,0,0,0)$ and $x^* = (55,55,55,55,27.5,27.5,55,55)$. The above UEP is not open-loop globally exponentially stable due to the existence of additional (congested) equilibria. This is shown in Figure 6, where the solution of the open-loop system, with constant inflows $v^* = (25,0,0,0,12.5,0,0,0)$, constant $d(t) \equiv (1,0,0,1,0.5)$ and $x_0 = [a_1,...,a_8]$, is attracted by the congested equilibrium $(111.8,111.8,111.8,111.8,27.5,27.5,92.82,92.82)'$ (Figure 6(a)) leading to outflow, which is 7.4 [veh] lower than the capacity flow of the 4$^{th}$ cell and 4.9 [veh] lower than the capacity flow of the 8$^{th}$ cell and a constant deviation of 125.5 [veh] for the Euclidean norm (Figure 6(b)). Therefore, if the objective is the operation of the freeway with largest possible outflow, then a control strategy will be needed.

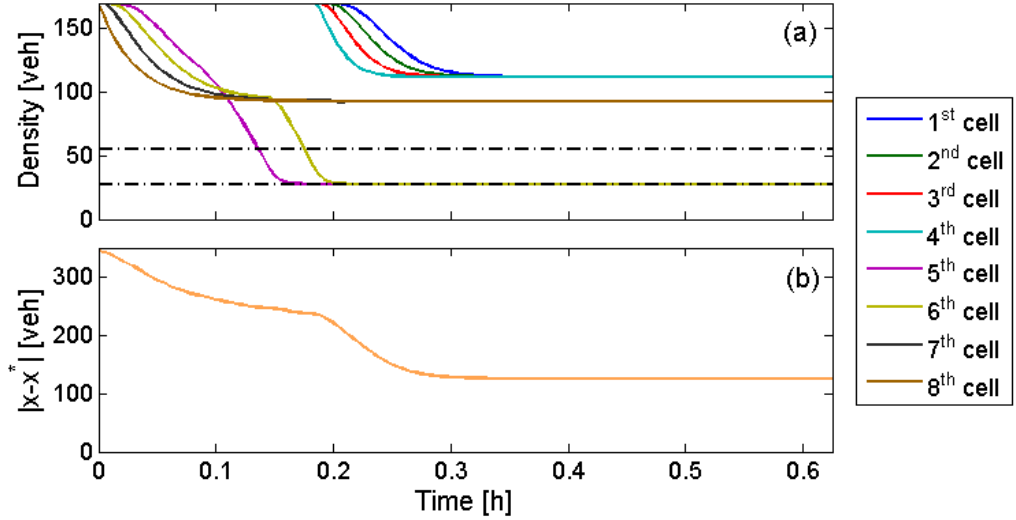

**Figure 6:** (a) The response of the density of every cell and (b) the evolution of the Euclidean norm of the deviation $x(t) - x^*$ of the state from the UEP, that is $|x(t) - x^*|$, for the open-loop system (5.2), (5.3), (5.4), (3.3), for initial condition $x_0 = [a_1,...,a_8]$, $v = v^*$ and $d(t) \equiv (1,0,1,0.5)$.

We constructed the matrix $K$ and the constant $\tau$ using the sufficient conditions provided from the proofs of the technical lemmas and propositions. Here, we simply used $K = 0.016 \cdot 1_{n \times n}$ and $\tau = 1/2$ which satisfy those conditions and allow for a good control performance with respect to overshooting effects. Figure 7 shows the response of the density of every cell for the closed-loop system (5.2), (5.3), (3.3) and three different initial conditions for constant $d(t) \equiv (1,0,0,0.5)$; Figure 7(a) is with $x_0 = (20,25,20,25,20,25,20,25)'$ corresponding to very low densities; Figure 7(b) is with $x_0 = (50,50,50,50,27,27,80,60)'$ corresponding to a more realistic traffic situation for which a sudden incident created congestion in a small part of the second freeway; and Figure 7(c) is with $x_0 = (a_1,...,a_8)'$ corresponding to a fully congested network. The feedback regulator is seen to respond very satisfactorily in these tests, exhibiting a fast convergence to the UEP for each one of the initial conditions.

Figure 8 shows again the response of the density of every cell for the closed-loop system (5.2), (5.3), (3.3) and three different initial conditions (same as those of Figure 7) for time-varying $d(t) = (d_1(t),...,d_4(t)) \in D$ taking values from a uniform distribution within $D = [0,1]^3 \times [0.22,0.3]$. In this case, although small oscillations exist, the rate of convergence to the UEP is similar to the previous test. This demonstrates the robustness of the feedback regulator (3.3) with respect to the uncertainties derived from the fundamental diagram (5.5), (5.6), (5.7).

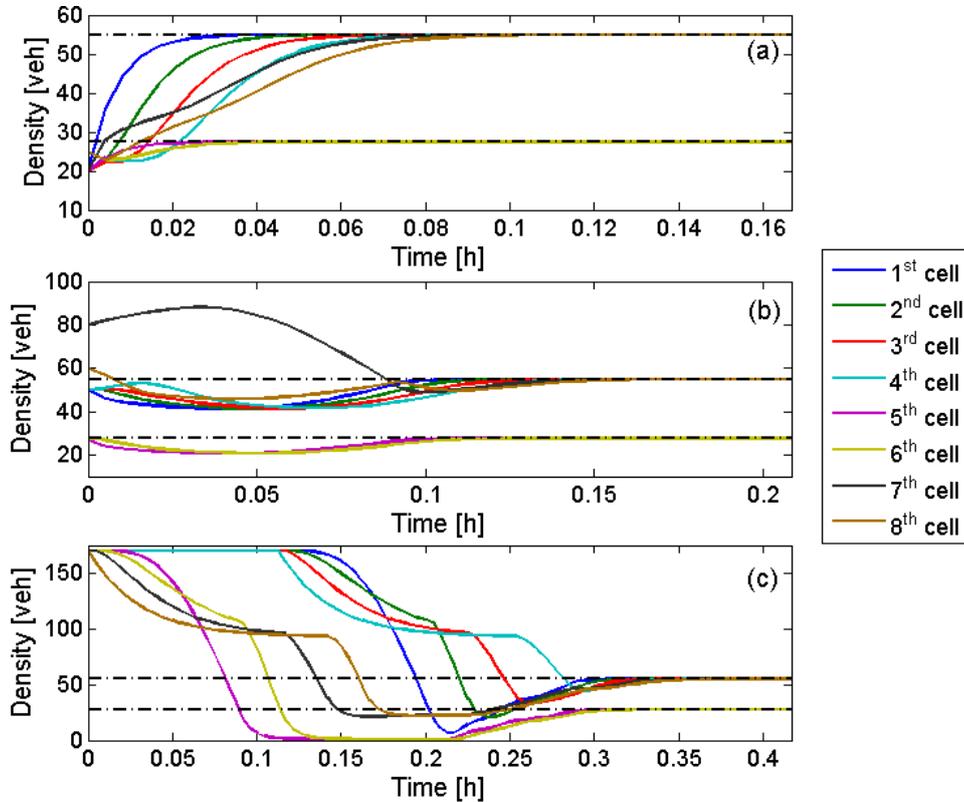

**Figure 7:** The response of the density of every cell for the closed-loop system (5.2), (5.3), (5.4), under the proposed feedback regulator (3.3) for different initial conditions; (a) $x_0 = (20,25,20,25,20,25,20,25)'$, (b) $x_0 = (50,50,50,50,27,27,80,60)'$ and (c) $x_0 = (a_1,...,a_8)'$ and for $d(t) \equiv (1,0,1,0.5)$.

Finally, Figure 9 shows the evolution of the Euclidean norm of the deviation of the state from the UEP, that is $|x(t) - x^*|$, for the closed-loop system (5.2), (5.3), (3.3) with $K = 0,016 \cdot 1_{n \times n}$ and $\tau = 1/2$, for different initial conditions; (a) $x_0 = (a_1,...,a_8)'$, (b) $x_0 = (150,140,60,120,120,100,160,130)'$, (c) $x_0 = (100,120,10,20,110,80,5,90)'$ and (d) $x_0 = (50,50,50,50,27,27,80,60)'$. By observing the evolution of norms, it can be concluded that the UEP, regardless of the initial conditions (as also guaranteed by the theoretical results), is reached within a small transient period.

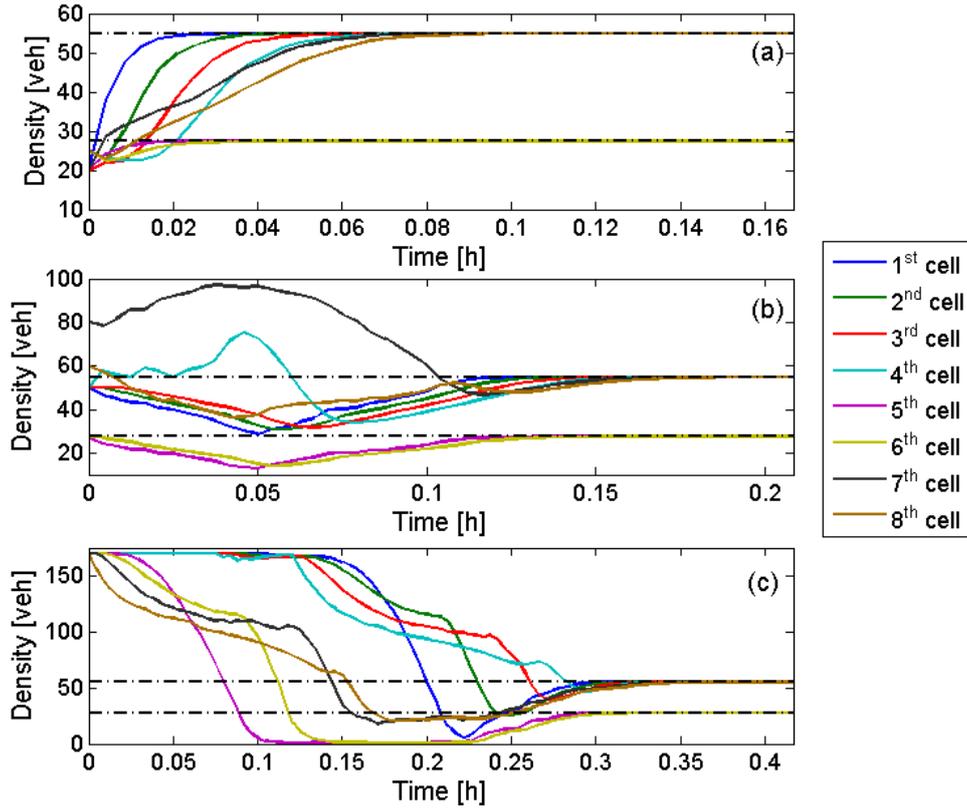

**Figure 8:** The response of the density of every cell for the closed-loop system (5.2), (5.3), (5.4), under the proposed feedback regulator (3.3) for different initial conditions; (a) $x_0 = (20,25,20,25,20,25,20,25)'$, (b) $x_0 = (50,50,50,50,27,27,80,60)'$ and (c) $x_0 = (a_1,..., a_8)'$ and for time-varying $d(t)$.

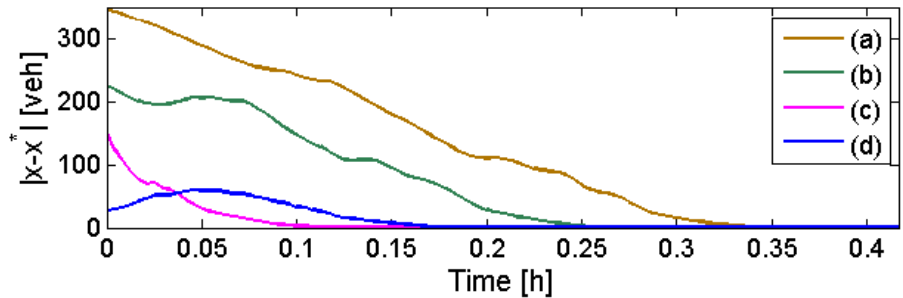

**Figure 9:** The time evolution of the Euclidean norm of the deviation $x(t) - x^*$ of the state from the UEP, that is $|x(t) - x^*|$, for the closed-loop system (5.2), (5.3), (5.4), under the proposed feedback regulator (3.3) and for time-varying $d(t)$ for different initial conditions; (a) $x_0 = (a_1,..., a_8)'$, (b) $x_0 = (150,140,60,120,120,100,160,130)'$, (c) $x_0 = (100,120,10,20,110,80,5,90)'$ and (d) $x_0 = (50,50,50,50,27,27,80,60)'$.

## 6. Concluding Remarks

This work provided a rigorous methodology for the construction of a parameterized family of explicit feedback laws that guarantee the robust global exponential stability of the UEP for general nonlinear uncertain discrete-time acyclic traffic networks. The construction of the global exponential feedback stabilizer is based on a vector Lyapunov function approach as well as certain important properties of acyclic traffic networks. The applicability of the obtained results to real control problems is demonstrated by conducting a simulation study, using a freeway-to-freeway network, with respect to various initial conditions. Simulation results demonstrate the efficacy of the proposed feedback control law with respect to the fast convergence to the UEP.

Future research will address robustness issues in a more comprehensive way. Moreover, future work includes application of the proposed methodology into an adaptive control framework, similarly to the work in [31], [32]. Moreover, testing the proposed feedback approach with other, more realistic (e.g. second-order) traffic simulation models, such as METANET [33] is also under way.

## Acknowledgments

The research leading to these results has received funding from the European Research Council under the European Union's Seventh Framework Programme (FP/2007-2013) / ERC Grant Agreement n. [321132], project TRAMAN21.

# Appendix

**Proof of Lemma 2.1:** Due to the fact that $P$ is strictly upper triangular, it holds that $p_{i,j} = 0$ for every $j \leq i$ and $p_{i,j} \in [0,1]$ otherwise. Therefore, in order to prove inequalities (2.10), it suffices to show that for every $r_n > 0$ there exist $r_i > 0$ ($i = 1,...,n-1$), such that $r_i > \sum_{j=i+1}^{n} r_j$ for $i = 1,...,n-1$. By choosing $r_i = 2^{n-i}$ we obtain:

$$2^{n-i} > \sum_{j=i+1}^{n} 2^{n-j} = \sum_{j=0}^{n-(i+1)} 2^j = 2^{n-i} - 1,$$

which holds for every $i = 1,...,n-1$. The proof is complete. ◁

**Proof of Lemma 2.2:** Due to the fact that $P$ is strictly upper triangular, it holds that $p_{i,j} = 0$ for every $j \leq i$ and $p_{i,j} \in [0,1]$ otherwise. Therefore, in order to prove inequalities (2.11), it suffices to show that there exist $\xi_i > 0$ ($i = 1,...,n$), such that the inequalities $\sum_{j=1}^{i-1} G_j \xi_j < L_i \xi_i$ hold for $i = 2,...,n$. For arbitrary $\xi_1 > 0$, we generate recursively the constants $\xi_i > 0$, $i = 2,...,n$, by using the following formula:

$$\xi_i = \frac{2}{L_i} \sum_{j=1}^{i-1} G_j \xi_j. \tag{A.0}$$

The proof is complete. ◁

**Proof of Lemma 2.3:** We prove this lemma by using the following Fact.

**Fact:** *The product of a strictly lower triangular matrix $A$ and a diagonal matrix $B$ is a strictly lower triangular matrix $C = AB$.*

**Proof of the Fact:** From the definition of matrix product, we have that $c_{i,j} = \sum_{k=1}^{n} a_{i,k} b_{k,j}$ for every $i, j \in \{1,...,n\}$. Due to the fact that $B$ is diagonal, it follows that $c_{i,j} = a_{i,j} b_{j,j}$. Due to the fact that $A$ is strictly lower triangular, it follows that $a_{i,j} = 0$ for $i \leq j$, which implies $c_{i,j} = 0$ for $i \leq j$. The proof of the Fact is complete. ◁

Since the matrix $P$ is a strictly upper triangular matrix, it holds that $P'$ is a strictly lower triangular matrix. Then, the product $P' diag(G)$ is a strictly lower triangular matrix from the above fact. Moreover, the matrix $I - diag(L)$ is a diagonal matrix. Then, the matrix $I + P' diag(G) - diag(L)$, which is the sum of a strictly lower triangular matrix and a diagonal matrix, is a lower triangular matrix with its diagonal elements being $(1-L_1,...,1-L_n)$ corresponding to its eigenvalues. Due to the fact that $L_i \in (0,1)$, we have that $\rho(I + P'diag(G) - diag(L)) = \max_i (|1-L_i|) = \max_i (1-L_i) < 1$. The proof is complete.  ◁

**Proof of Proposition 2.4:** Using (2.9) we obtain

$$\left(\sum_{i=1}^{n} r_i x_i\right)^+ = \sum_{i=1}^{n} r_i x_i + \sum_{i=1}^{n} r_i w_i(d,x,v) v_i - \sum_{i=1}^{n} r_i s_i(d,x,v) f_i(d,x_i) + \sum_{i=1}^{n}\sum_{j=1}^{n} r_i p_{j,i} s_j(d,x,v) f_j(d,x_j)$$

$$= \sum_{i=1}^{n} r_i x_i + \sum_{i=1}^{n} r_i w_i(d,x,v) v_i - \sum_{i=1}^{n}\left(1 - \sum_{j=1}^{n} r_i^{-1} r_j p_{i,j}\right) r_i s_i(d,x,v) f_i(d,x_i).$$
(A.1)

Define

$$Q := \min_{i=1,...,n}\left(1 - \sum_{j=1}^{n} r_i^{-1} r_j p_{i,j}\right).$$
(A.2)

Notice that Lemma (2.1) guarantees that $Q > 0$. Using (A.1), (A.2) and the fact that $w_i(d,x,v) \in [0,1]$, we obtain:

$$\left(\sum_{i=1}^{n} r_i x_i\right)^+ \leq \sum_{i=1}^{n} r_i x_i + \sum_{i=1}^{n} r_i v_i - Q\sum_{i=1}^{n} r_i s_i(d,x,v) f_i(d,x_i).$$
(A.3)

Assumption (H1) guarantees that $f_i(d,x_i) \geq L_i x_i$ for all $x_i \in [0,\tilde{\delta}_i]$ and $i=1,...,n$. Moreover, assumption (H1) in conjunction with the previous inequality guarantees that $f_i(d,x_i) \geq \min\left(f_i^{\min}, L_i \tilde{\delta}_i\right)$ for all $x_i \in [\tilde{\delta}_i, a_i]$ and $i=1,...,n$. Combining we obtain $f_i(d,x_i) \geq \theta_i x_i$ for all $x_i \in [0,a_i]$ and $i=1,...,n$, where $\theta_i := \min\left(L_i, a_i^{-1} f_i^{\min}, a_i^{-1} L_i \tilde{\delta}_i\right)$. Notice that $\theta_i > 0$ for $i=1,...,n$ and define $\Theta := \min_{i=1,...,n}(\theta_i)$. It follows that

$$f_i(d,x_i) \geq \Theta x_i \text{ for all } x_i \in [0,a_i] \text{ and } i=1,...,n.$$
(A.4)

Again, notice that $\Theta > 0$. We obtain from (A.3) and (A.4):

$$\left(\sum_{i=1}^{n} r_i x_i\right)^+ \leq \sum_{i=1}^{n} r_i x_i + \sum_{i=1}^{n} r_i v_i - Q\Theta \sum_{i=1}^{n} r_i s_i(d,x,v) x_i.$$
(A.5)

Since the set $D \times S$ is compact, it follows from continuity of the functions $g_i : D \times S \to \Re_+$ ($i=1,...,n$) that there exists $\delta > 0$ such that $|g_i(d,x) - g_i(\tilde{d},0)| < \frac{1}{2}\min_{i=1,...,n}(\tilde{\varepsilon}_i)$ for all $i=1,...,n$, $(d,x) \in D \times S$, $\tilde{d} \in D$ with $|d-\tilde{d}| + |x| < \delta$. Since $p_{i,j} \in [0,1]$, for $i,j = 1,...,n$, $v_i \leq \min\{g_i(d,0) : d \in D\} - \tilde{\varepsilon}_i$ for $i=1,...,n$, and since $0 < f_i(d,z) < z$ for all $z \in (0,a_i]$ and $i=1,...,n$ (see Assumption (H1)), we get for all $i=1,...,n$, $(d,x) \in D \times S$ with $|x| < \min\left(\delta, \frac{1}{2n}\min_{i=1,...,n}(\tilde{\varepsilon}_i)\right)$:

$$v_i + \sum_{j=1}^{n} p_{j,i} f_j(d,x_j) \leq \min\{g_i(d,0) : d \in D\} - \tilde{\varepsilon}_i + \sum_{j=1}^{n} x_j \leq \min\{g_i(d,0) : d \in D\} - \frac{\tilde{\varepsilon}_i}{2} \leq g_i(d,x).$$

It follows from the above inequality, (2.8) and (A.5) that the following inequality holds for all $(d,x) \in D \times S$ with $|x| < \min\left(\delta, \frac{1}{2n} \min_{i=1,\ldots,n}(\widetilde{\varepsilon}_i)\right)$:

$$\left(\sum_{i=1}^{n} r_i x_i\right)^+ \leq (1-Q\Theta)\sum_{i=1}^{n} r_i x_i + \sum_{i=1}^{n} r_i v_i. \qquad (A.6)$$

Define $V = \{v = (v_1,\ldots,v_n)' \in \Re_+^n : v_i \leq \min(v_i^{\max}, \min\{g_i(d,0) : d \in D\}) - \widetilde{\varepsilon}_i, i=1,\ldots,n\}$. We next claim that there exists a constant $\gamma > 0$ such that

$$\sum_{i=1}^{n} r_i \widetilde{s}_i(d,x,v) x_i \geq \gamma \sum_{i=1}^{n} r_i x_i, \text{ for all } (d,x,v) \in D \times S \times V \text{ with } |x| \geq \min\left(\delta, \frac{1}{2n} \min_{i=1,\ldots,n}(\widetilde{\varepsilon}_i)\right). \qquad (A.7)$$

Indeed, we define

$$\gamma := \inf\left\{\frac{\sum_{i=1}^{n} r_i \widetilde{s}_i(d,x,v) x_i}{\sum_{i=1}^{n} r_i x_i} : (d,x,v) \in D \times S \times V, |x| \geq \min\left(\delta, \frac{1}{2n}\min_{i=1,\ldots,n}(\widetilde{\varepsilon}_i)\right)\right\}. \qquad (A.8)$$

In order to show (A.7), it suffices to show that $\gamma \geq 0$ as defined by (A.8) is positive. Continuity of the functions $\widetilde{s}_i : D \times S \times \Re_+^n \to [0,1]$ ($i=1,\ldots,n$) and compactness of the set $\left\{(d,x,v) \in D \times S \times V : |x| \geq \min\left(\delta, \frac{1}{2n}\min_{i=,\ldots,n}(\widetilde{\varepsilon}_i)\right)\right\}$ imply that there exists $(d^*, x^*, v^*) \in D \times S \times V$ with $|x^*| \geq \min\left(\delta, \frac{1}{2n}\min_{i=1,\ldots,n}(\widetilde{\varepsilon}_i)\right)$ and

$$\sum_{i=1}^{n} r_i \widetilde{s}_i(d^*, x^*, v^*) x_i^* = \gamma \sum_{i=1}^{n} r_i x_i^*. \qquad (A.9)$$

We proceed by using a contradiction argument. Suppose that $\gamma = 0$. It follows from (A.9) that $x_i^* \widetilde{s}_i(d^*, x^*, v^*) = 0$ for $i=1,\ldots,n$. However, (H3) and (2.13) imply that $x^* = 0$, which contradicts the fact that $|x^*| \geq \min\left(\delta, \frac{1}{2n}\min_{i=1,\ldots,n}(\widetilde{\varepsilon}_i)\right)$.

It follows from (A.5), (A.6), (A.7) and the fact that $s_i(d,x,v) \geq \widetilde{s}_i(d,x,v)$ for all $(d,x,v) \in D \times S \times \Re_+^n$, $i=1,\ldots,n$, that inequality (2.14) holds with $C := Q\Theta \min(1,\gamma)$. The proof is complete. ◁

**Proof of Lemma 3.4:** Proposition 3.3 guarantees that there exist constants $\beta_i \in (x_i^*, \mu_i]$ ($i=1,\ldots,n$) such that (3.5) holds. Using (4.3), (4.5), (3.1), we obtain for all $x \in \Omega$ and $i=1,\ldots,n$:

$$\begin{aligned}x_i^+ - x_i^* &\leq -(v_i^* - b_i)\min\left(1, \tau^{-1}\sum_{j=1}^{n} K_{i,j}\max(0, x_j - x_j^*)\right) \\ &+ (1-L_i)\max(0, x_i - x_i^*) + \sum_{j=1}^{n} p_{j,i} G_j \max(0, x_j - x_j^*)\end{aligned} \qquad (A.10)$$

Inequality (3.6) is a direct consequence of inequalities (A.10). Using Assumption (H1), we get $G_i(x_i - x_i^*) \leq f_i(d, x_i) - f_i(d, x_i^*) \leq L_i(x_i - x_i^*)$ for $i=1,\ldots,n$ and $x_i \leq x_i^*$. Notice that assumption (H1) and the fact that $\mu_i \in (0, \widetilde{\delta}_i)$ ($i=1,\ldots,n$) guarantee that the mappings $x_i \to x_i - f_i(d, x_i)$, $x_i \to f_i(d, x_i)$ are non-decreasing for $x_i \in [0, \beta_i]$, $i=1,\ldots,n$. It follows that:

$$f_i(d, x_i) - f_i(d, x_i^*) \geq -G_i \max(0, x_i^* - x_i), \quad x_i - x_i^* + f_i(d, x_i^*) - f_i(d, x_i) \geq -(1 - L_i) \max(0, x_i^* - x_i)$$
$$\text{for } x_i \in [0, \beta_i], \; i = 1, \ldots, n. \tag{A.11}$$

Using (4.3), (3.1), we obtain for all $x \in \Omega$, $i = 1, \ldots, n$:

$$x_i^+ - x_i^* \geq -(1 - L_i) \max(0, x_i^* - x_i) - (v_i^* - b_i) \min\left(1, \tau^{-1} \sum_{j=1}^n K_{i,j} \max(0, x_j - x_j^*)\right)$$
$$- \sum_{j=1}^n p_{j,i} G_j \max(0, x_j^* - x_j) \tag{A.12}$$

Inequalities (A.12) imply the following inequalities for all $x \in \Omega$ and $i = 1, \ldots, n$:

$$x_i^* - x_i^+ \leq (1 - L_i) \max(0, x_i^* - x_i) + (v_i^* - b_i) \min\left(1, \tau^{-1} \sum_{j=1}^n K_{i,j} \max(0, x_j - x_j^*)\right)$$
$$+ \sum_{j=1}^n p_{j,i} G_j \max(0, x_j^* - x_j) \tag{A.13}$$

Inequality (3.7) is a direct consequence of inequalities (A.13) and the fact that $\min\left(1, \tau^{-1} \sum_{j=1}^n K_{i,j} \max(0, x_j - x_j^*)\right) \leq \tau^{-1} \sum_{j=1}^n K_{i,j} \max(0, x_j - x_j^*)$ for all $x \in \Omega$ and $i = 1, \ldots, n$. The proof is complete. ◁

**Proof of Lemma 3.5:** Lemma 3.4 guarantees that there exist constants $\beta_i \in (x_i^*, \mu_i]$ ($i = 1, \ldots, n$) such that (3.5), (3.6), (3.7) hold. Using (3.6), (3.7) in conjunction with the following inequalities which hold for all $x \in \Re^n$

$$1'_n h(x - x^*) + 1'_n h(x^* - x) = \sum_{i=1}^n |x_i - x_i^*| \geq |x - x^*| \tag{A.14}$$

$$h(x - x^*) \leq |x - x^*| 1_n, \quad h(x^* - x) \leq |x - x^*| 1_n \tag{A.15}$$

and the facts that $(I + P' \text{diag}(G) - \text{diag}(L))$, $K \in \Re_+^{n \times n}$, $\text{diag}(v^* - b)$ are non-negative matrices and $\tau > 0$, we are in a position to guarantee the existence of $\delta > 0$ such that

$$|x^+ - x^*| \leq 1'_n \left(2I + 2P' \text{diag}(G) - 2\text{diag}(L) + \text{diag}(v^* - b) \tau^{-1} K\right) 1_n |x - x^*|,$$
for all $x \in S, d \in D$ with $|x - x^*| < \delta$ and $v = v^* - \text{diag}(v^* - b)(1_n - h(1_n - \tau^{-1} K h(x - x^*)))$. $\tag{A.16}$

Since $x^+ \in S$ for $x \in S, d \in D$, $v = (v_1, \ldots, v_n)' \in \Re_+^n$ and since $S = [0, a_1] \times \ldots \times [0, a_n]$, it follows that:

$$|x^+ - x^*| \leq 2\sqrt{n} \max_{i=1,\ldots,n}(a_i),$$
for all $x \in S, d \in D$ and $v = v^* - \text{diag}(v^* - b)(1_n - h(1_n - \tau^{-1} K h(x - x^*)))$. $\tag{A.17}$

Estimates (A.16), (A.17) imply that (3.8) holds with $M := 1'_n \left(2I + 2P' \text{diag}(G) - 2\text{diag}(L) + \text{diag}(v^* - b) \tau^{-1} K\right) 1_n + 2\sqrt{n} \delta^{-1} \max_{i=1,\ldots,n}(a_i)$. The proof is complete. ◁

**Proof of Lemma 3.6:** Lemma 3.4 guarantees that there exist constants $\beta_i \in (x_i^*, \mu_i]$ ($i=1,...,n$) such that (3.5), (3.6), (3.7) hold. Let $K \in \Re_+^{n \times n}$ be a matrix so that

$$K_{i,j} \geq \frac{1}{\min_k (\beta_k - x_k^*)}, \text{ for every } i,j = 1,...,n. \tag{A.18}$$

It follows from (A.18) that for every $\tau \in (0,1)$ we have

$$\tau^{-1} K_{i,j} (\beta_j - x_j^*) \geq 1, \text{ for every } i,j = 1,...,n. \tag{A.19}$$

Using the fact that $x_i^* < \beta_i$ for every $i=1,...,n$, it follows from (A.19) and the fact that $\max(0, x_i - x_i^*) = x_i - x_i^*$ for $x \in S \setminus \Omega$ (recall that $\Omega = [0, \beta_1] \times \cdots \times [0, \beta_n]$), that

$$\tau^{-1} \sum_{j=1}^{n} K_{i,j} \max(0, x_j - x_j^*) \geq 1, \text{ for every } i=1,...,n \text{ and } x \in S \setminus \Omega \tag{A.20}$$

and consequently, since $v = v^* - diag(v^* - b)(1_n - h(1_n - \tau^{-1} Kh(x - x^*)))$, we get

$$v = b, \text{ for } x \in S \setminus \Omega. \tag{A.21}$$

In order to show that the set $\Omega = [0, \beta_1] \times ... \times [0, \beta_n]$ is a TR for the closed-loop system (2.9) with (3.3), it suffices to show that for every $x_0 \in S$, $\{d(t) \in D\}_{t=0}^{\infty}$ the solution $x(t)$ of the closed-loop system (2.9) with (3.3) and initial condition $x(0) = x_0$ corresponding to input $\{d(t) \in D\}_{t=0}^{\infty}$ satisfies $x(t) \in \Omega$ for all $t \geq m$, where

$$m := \left\lceil \frac{\ln\left(C \min_{i=1,...,n}(r_i \beta_i) - r'b\right) - \ln(Cr'a)}{\ln(1-C)} \right\rceil + 1 \tag{A.22}$$

where $a = (a_1,...,a_n)' \in \text{int}(\Re_+^n)$. We proceed by contradiction. Suppose that there exist $x_0 \in S$, $\{d(t) \in D\}_{t=0}^{\infty}$ such that the solution $x(t)$ of the closed-loop system (2.9) with (3.3) and initial condition $x(0) = x_0$ corresponding to input $\{d(t) \in D\}_{t=0}^{\infty}$ satisfies $x(t) \notin \Omega$ for certain $t \geq m$. Since the set $\Omega = [0, \beta_1] \times ... \times [0, \beta_n]$ is positively invariant (a direct consequence of (3.5)), it follows that $x(q) \notin \Omega$ for all $q = 0,1,...,m$. Define

$$I(q) := r'x(q) \tag{A.23}$$

and notice that (2.14), (A.21) imply the following estimate for all $q = 0,1,...,m$:

$$I(q+1) \leq (1-C)I(q) + r'b. \tag{A.24}$$

Estimate (A.24) implies the following estimate for all $q = 0,1,...,m+1$:

$$I(q) \leq (1-C)^q I(0) + C^{-1} r'b \left(1 - (1-C)^q\right). \tag{A.25}$$

Since $I(0) = r'x(0) = r'x_0 \leq r'a$ for all $x_0 \in S$, we obtain from (A.25) for all $q = 0,1,...,m+1$:

$$I(q) \leq (1-C)^q r'a + C^{-1} r'b. \tag{A.26}$$

Estimate (A.26) in conjunction with definition (A.22) implies that $I(m) \leq \min_{i=1,...,n}(r_i \beta_i)$, which combined with definition (A.23) shows that $x(m) \in \Omega$, a contradiction. The proof is complete. ◁